\documentclass[11pt,amsfonts]{article}
\usepackage{graphicx}
\usepackage{latexsym}
\usepackage{amssymb}
\usepackage{amsmath}
\usepackage{enumerate}
\usepackage{color}
\usepackage{layout}

\newtheorem{prop}{Proposition}
\newtheorem{lemma}{Lemma}

\newtheorem{theorem}{Theorem}
\newtheorem{remark}{Remark}

\def\real{{\mathord{{\rm I\kern-2.8pt R}}}}        
\def\inte{{\mathord{{\rm I\kern-2.8pt N}}}}

\def\sZZ{{\rm Z\kern-2.8ptem{}Z}}

\def\z{{\mathchoice
		{\sZZ}
		{\sZZ}
		{\rm Z\kern-0.30em{}Z}
		{\rm Z\kern-0.25em{}Z} }}
\def\sQQ{{\kern 0.27em \vrule height1.45ex width0.03em depth0em
		\kern-0.30em \rm Q}}
\def\qu{{\mathchoice
		{\sQQ}
		{\sQQ}
		{\kern 0.225em \vrule height1.05ex width0.025em depth0em \kern-0.25em \rm Q}
		{\kern 0.180em \vrule height0.78ex width0.020em depth0em \kern-0.20em \rm Q}
}}
\def\sCC{{\kern 0.27em \vrule height1.45ex width0.03em depth0em
		\kern-0.30em \rm C}}
\def\complex{{\mathchoice
		{\sCC}
		{\sCC}
		{\kern 0.225em \vrule height1.05ex width0.025em depth0em \kern-0.25em \rm C}
		{\kern 0.180em \vrule height0.78ex width0.020em depth0em \kern-0.20em \rm C}
}}

\newcommand{\ba}{\begin{array}}
	\newcommand{\ea}{\end{array}}
\newcommand{\be}{\begin{equation}}
	\newcommand{\ee}{\end{equation}}
\newcommand{\bea}{\begin{eqnarray}}
	\newcommand{\eea}{\end{eqnarray}}
\newcommand{\beaa}{\begin{eqnarray*}}
	\newcommand{\eeaa}{\end{eqnarray*}}

\def\z{\zeta}

\font\tenmath=msbm10 \font\sevenmath=msbm7 \font\fivemath=msbm5
\newfam\mathfam \textfont\mathfam=\tenmath
\scriptfont\mathfam=\sevenmath \scriptscriptfont\mathfam=\fivemath

\def \={{\buildrel {\rm (law)} \over =}}

\def\qed{ \hfill \vrule width.25cm height.25cm depth0cm\smallskip}

\newcommand{\basa}{\begin{assumption}}
	\newcommand{\easa}{\end{assumption}}

\newcommand{\bas}{\begin{assum}}
	\newcommand{\eas}{\end{assum}}

\newcommand{\ignore}[1]{}
\begin{document}
	
	\renewcommand{\thefootnote}{\fnsymbol{footnote}}
	
	\renewcommand{\thefootnote}{\fnsymbol{footnote}}

	\title{ Temporal quadratic and higher order  variation for  the nonlinear stochastic heat equation and applications to parameter estimation}
	\author{Christian Olivera \footnote{
			C. Olivera is partially supported by FAPESP by the grant  $2020/04426-6$,  
			by FAPESP-ANR by the grant Stochastic and Deterministic Analysis for Irregular Models$-2022/03379-0$ 
			and CNPq by the grant $422145/2023-8$. } $^{1,}$ \hskip0.2cm 
		Ciprian A. Tudor \footnote{	 C. Tudor also acknowledges support from   the  ANR project SDAIM 22-CE40-0015,  ECOS SUD (project C2107), Japan Science and Technology Agency CREST  (grant JPMJCR2115) and  by the Ministry of Research, Innovation and Digitalization (Romania), grant CF-194-PNRR-III-C9-2023.}$^{2,3}$ \vspace*{0.1in} \\
		$^{1}$ Departamento de Matem\'atica, Universidade Estadual de Campinas,\\
		13.081-970-Campinas-SP-Brazil. \\
		colivera@ime.unicamp.br \vspace*{0.1in} \\
		$^{2}$  Universit\'e de Lille, CNRS\\
		Laboratoire Paul Painlev\'e UMR 8524\\
		F-59655 Villeneuve d'Ascq, France.\\
		$^{3}$ Bucharest University for Economic Studies, Romania.\\
		\quad tudor@math.univ-lille1.fr\vspace*{0.1in}}
	\maketitle
	
	\maketitle
	
	\begin{abstract}
		We consider the stochastic heat equation  which includes a fractional power of the Laplacian of order $\alpha \in (1, 2]$ and it is driven by a nonlinear space-time Gaussian white noise. We study two types of power variations for the solution to this equation: the renormalized quadratic variation and the power variation of order $\frac{2\alpha}{\alpha -1}$, both over an equidistant partition of the unit interval. We prove that these two sequences admit nontrivial limits when the mesh of the partition goes to zero.  We apply these results to identify certain  parameters of the stochastic heat equation.   
	\end{abstract}
	\vskip0.3cm
	
	{\bf 2010 AMS Classification Numbers:}  60G15, 60H05, 60G18.

	\vskip0.3cm

	{\bf Key Words and Phrases}: Stochastic heat equation, Fractional Brownian motion, Fractional Laplacian, power variation;    parameter estimation.

	\section{Introduction}
	Our work concerns the  following  stochastic heat equation
	\begin{equation}
		\label{intro-1}
		\frac{ \partial u_{\theta}}{\partial t} (t,x)= -\theta(-\Delta) ^{\frac{ \alpha }{2}} u_{\theta }(t,x)+  \sigma (u_{\theta}(t,x)) \dot{W}(t,x), \hskip0.4cm t\geq 0, x\in \mathbb{R},
	\end{equation}
	with vanishing initial value $ u_{\theta }(0,x)=0$ for every $x\in \mathbb{R}$. The random perturbation $W$ in (\ref{intro-1}) is a space-time white noise, $ -(-\Delta) ^{\frac{ \alpha }{2}}$ stands for the fractional Laplacian operator of order $\alpha \in (1, 2]$ and $ \sigma $ is a nonlinear Lipschitz function. The purpose is to study the limit behavior of the power variations (with respect to the time variable) of the solution to (\ref{intro-1}) and to apply these results  to estimate the  two  parameters involved in the model (\ref{intro-1}), i.e. the drift parameter $\theta>0$ and the index $\alpha $ of the fractional Laplacian. We will assume that the solution $ u_{\theta}$ to (\ref{intro-1}) is observed at discrete times on the unit interval $[0,1]$ and at a fixed spatial point, i.e. we dispose on the observations $ (u_{\theta}(t_{i}, x), i=0,1,...,N)$, where $ t_{i}=\frac{i}{N}$ for $i=0,1,...,N$ and $x\in \mathbb{R}$ is an arbitrary spatial point.  The parameter estimation for stochatic partial differential equations (SPDEs in the sequel), in particular   via power variations,  constitutes nowadays a very active research direction. The reader may consult the website https://sites.google.com/view/stats4spdes for a  vast bibliography on this topic.  
	
	Many of these works concern the situation of linear SPDEs, i.e. the random noise is additive and the nonlinear coefficient $\sigma$  is identically one. The reader may consult \cite{Cia} or the recent monograph \cite{T2}  for a survey of the estimation techniques for SPDEs with additive Gaussian noise. Other references on statistical inference for linear SPDEs are, among others, \cite{AT1}, \cite{BT1}, \cite{BT2}, \cite{CDK}, \cite{CK}, \cite{CKP}, \cite{G1}, \cite{GR},  \cite{GT}, \cite{HiTr}, \cite{Mahtu2}, \cite{KU}, \cite{TKU} or \cite{ZZ}.  The scientific literature for the statistical inference on SPDEs with nonlinear Gaussian noise is less rich. Among the relatively few works dealing with these aspects, we mention the work \cite{PoTr} for the study of $p$-variations of the solution to the nonlinear heat equation driven by a space-time white noise,  \cite{Chong} for the study of the temporal  $p$-variation of the nonlinear heat equation with standard Laplacian  driven by a Gaussian noise with spatial correlation, \cite{CD} for the analysis of the $p$-variation in time of the solution to the nonlinear fractional stochastic heat equation on a bounded spatial domain and \cite{GT} for the spatial $p$-variations of the solution to (\ref{intro-1}). See also \cite{CH} or \cite{GR}  for a statistical analysis of some nonlinear SPDEs.

	We will actually analyze two types of variations (in time)  for the solution of the SPDE (\ref{intro-1}): the renormalized temporal quadratic variation (denoted by $ V_{N,x}(u_{0})$ and given by (\ref{10f-1})) and the temporal $\frac{2\alpha}{\alpha-1}$-variation (the sequence denoted by $U_{N,x}(u_{0})$ and defined by (\ref{10f-2})). The analysis   of the second sequence (the "standard" power variation $U_{N,x}(u_{0})$) constitutes a natural problem and it has been studied in the literature for various SPDEs (see e.g. \cite{CD}, \cite{MahTu}, \cite{GT}, among many others). On the other hand, the "renormalized" quadratic variation $ V_{N, x}(u_{0})$, although less used in the literature, appears to be very useful in statistical inference. In particular, it allows to construct a consistent estimator for the parameter $\alpha$ appearing in the definition of the fractional Laplacian. 
	
	Based on the limit behavior of these two sequences, we are able to construct two consistent estimators for the drift   parameter $\theta $  and also a consistent estimator for the anomality parameter $\alpha$. We notice that the estimation of $\alpha$, although it seems a natural problem, has not have been treated in the literature, as far as we know. 
	
	The main idea of the proof of the limit theorem satisfied by  the variations of the solution is to approximate the increments of the solution to (\ref{intro-1}) by the increments of the solution to the corresponding linear heat equation (with $\sigma \equiv 1$) and then to use the fact that this later behaves as a perturbed fractional Brownian motion (fBm in the sequel), i.e. the sum of the fBm and  of  a Gaussian process with rather regular sample paths.  A similar idea has been used in \cite{GT} for the study of the spatial renormalized quadratic variation of the solution to the SPDE (\ref{intro-1}), but the perturbation of the fBm appearing in that context has different properties.

	We structured our work as follows.  In Section 2 we included the presentation of the fractional stochastic heat equation and of its solution. We also included some preliminary results on the quadratic variation of the fBm and of the perturbed fBm. In Section 3 we study the limit behavior of the quadratic and higher order variations of the solution to the stochastic heat equation with nonlinear multiplicative noise. Section 4 contains the applications to statistical inference while in Section 5 (the Appendix) we recall the basic properties of the Green kernel associated with the fractional Laplacian and we give the proof of some  technical results.

	\section{Preliminaries}
	In this preliminary part, we present some basic facts concerning the stochastic heat equation with fractional Laplacian  and its mild solution. Based on some known results for fBm, we also deduce some new properties  concerning the increments and the quadratic variation of a perturbed fBm. 
	
	\subsection{The solution}
	We start by describing the model (\ref{intro-1}) without drift parameter, i.e. we assume for the moment that $ \theta =1$. That is, 	we deal with the following SPDE
	\begin{equation}
		\label{2}
		\frac{\partial}{\partial t} u(t,x)= -(-\Delta ) ^{\frac{\alpha}{2}}u(t,x) + \sigma (u (t,x))\dot{W} (t,x), \hskip0.4cm t\geq 0, x\in \mathbb{R}.
	\end{equation}
	In (\ref{2}), the random perturbation $W$ is assumed to be  a space-time white noise, i.e. a centered Gaussian field $\left( W (t,A), t\geq 0, A \in \mathcal{B}_{b}(\mathbb{R})\right)$ with covariance function given by 
	\begin{equation}\label{covW}
		\mathbf{E} W(t,A) W(s, B) = (t\wedge s) \lambda (A \cap B), \hskip0.3cm s,t \geq 0, A, B \in \mathcal{B}_{b} (\mathbb{R}).
	\end{equation}
	We denoted by $\lambda$ the Lebesgue measure on $\mathbb{R}$ and $ \mathcal{B}_{b}(\mathbb{R})$  stands for the set of bounded Borel subsets of the real line $\mathbb{R}$. The operator $ -(-\Delta ) ^{\frac{\alpha}{2}}$ is the fractional Laplacian of order $\frac{\alpha}{2}$, where $\alpha \in (1, 2]$. It  reduces to  the standard Laplacian when $\alpha=2$.  We will mainly use the properties of the Green kernel associated to this operator which are recalled in the Appendix. The diffussion coefficient $\sigma$ in (\ref{2}) is assumed to be (globally) Lipschitz continuous, i.e. there exists $K>0$ such that, for every $x,y\in \mathbb{R}$,
	\begin{equation}\label{2i-3}
		\vert \sigma (x)- \sigma (y) \vert\leq K \vert x-y\vert.
	\end{equation}
	We know (see e.g. \cite{ANTV} or \cite{DD}) that, if (\ref{2i-3}) holds true, then the stochastic partial differential equation (\ref{2}) admits a unique mild solution given by 
	\begin{equation}
		\label{mild1}
		u(t,x)= \int_{0} ^{t} \int_{\mathbb{R}} G_{\alpha} (t-s, x-y)\sigma (u(s,y)) W (ds, dy), 
	\end{equation}
	where $G_{\alpha}$ stands for the Green kernel associated to the fractional heat equation  (see formula (\ref{ga}) in the Appendix). The stochastic integral with respect to the noise $W$ from the above formula (\ref{mild1}) is the so-called Dalang-Walsh integral (see \cite{Da} and \cite{Walsh}). Recall that the integral 
	\begin{equation*}
		\int_{0}^{\infty} \int_{\mathbb{R}} X(s,y) W(ds,dy)
	\end{equation*}
	is well-defined for random fields $( X(t,x), t\geq 0, x\in \mathbb{R})$ which are adapted to the filtration $\left( \mathcal{F}_{t}, t\geq 0\right), \mathcal{F}_{t}=\sigma\{ W(s,A), 0\leq s\leq t, A\in \mathcal{B}_{b}(\mathbb{R})\}$ and satisfy 
	\begin{equation*}
		\mathbf{E}	\int_{0} ^{\infty} \int_{\mathbb{R}} X (s,y) ^{2} dyds <\infty.
	\end{equation*}
	The Dalang-Walsh stochastic integral satisfies the following  It\^o-type  isometry 
	\begin{equation}\label{iso-dw}
		\mathbf{E}\left(\int_{0} ^{\infty} \int_{\mathbb{R}} X(s,y) W (ds, dy)\right) ^{2}= 	\mathbf{E} \int_{0} ^{\infty} \int_{\mathbb{R}} X (s,y) ^{2} dyds. 
	\end{equation}
	
	We will also need in the sequel some properties of the mild solution (\ref{mild1}). For their proofs, we refer to e.g. \cite{ANTV} or \cite{DD}. Let $T>0$.  Then, if $ C_{T} $ denotes a strictly positive constant depending on $T$, 
	\begin{enumerate}
		\item For every $p\geq 1$, we have
		\begin{equation}\label{3i-5}
			\sup_{t\in [0,T], x\in \mathbb{R}}		\mathbf{E}\left| u(t,x)\right| ^{p} \leq C_{T}.
		\end{equation}
		
		\item For every $s, t\in [0, T], x,y\in \mathbb{R}$ and for every $p\geq 2$,
		\begin{equation}\label{2i-2}
			\mathbf{E} \left| u(t,x)- u(s, y) \right| ^{p} \leq C_{T} \left[ \vert t-s\vert ^{ (1-\frac{1}{\alpha})\frac{p}{2}}+ \vert x-y\vert ^{(\alpha-1) \frac{p}{2}} \right].
		\end{equation}
		In particular, the sample paths of $u$ are H\"older continuous of order $\delta \in \left( 0, \frac{ \alpha-1}{2\alpha }\right) $ in time and of order $ \delta \in \left(0, \frac{ \alpha-1}{2}\right)$ in space. 
	\end{enumerate}

	\subsection{Some useful properties of the fractional Brownian motion}
	Since the solution to the (linear) stochastic heat equation is related to the fractional Brownian motion, let us start with some results concerning fBm. Recall that the fBm with Hurst  parameter $H\in (0, 1)$ is defined as a centered Gaussian process $(B^{H}_{t}, t\geq 0)$ with covariance function 
	\begin{equation*}
		\mathbf{E} B ^{H}_{t} B ^{H}_{s} =\frac{1}{2} \left( t^{2H}+ s^{2H}-\vert t-s\vert ^{2H}\right), \mbox{ for every } s, t\geq 0.
	\end{equation*}
	
	In the sequel, let
	\begin{equation}
		\label{ti}t_{i}= \frac{i}{N}, \mbox{ for } i=0, 1,...,N.
	\end{equation}

	\begin{lemma}\label{ll11}
		Let $(B ^{H}_{t}, t\geq 0)$ be a fBm with Hurst index $H\in (0,1)$. Then, for every $i=0,\ldots, N-1$, 
		\begin{equation*}
			\mathbf{E} \left[ N ^{2H-1} \left( B ^{H} _{t_{i+1}} - B ^{H} _{t_{i}} \right)^{2} -\frac{1}{N} \right] ^{2}= \frac{2}{N^{2}}.
		\end{equation*}
		
	\end{lemma}
	{ \bf Proof: } This result has been obtained in Lemma 3 in \cite{GT}. \qed

	\begin{prop}\label{pp2}
		Let $ (B ^{H}_{t}, t\geq 0)$ be a fractional Brownian motion with Hurst index $H\in \left(0, \frac{3}{4}\right)$. For $N\geq 1$, denote
		\begin{equation*}
			V_{N}(B ^{H})= N ^{2H-1} \sum_{i=0} ^{N-1} \left( B^{H}_{ t_{i+1}} - B ^{H}_{t_{i}} \right) ^{2}.
		\end{equation*}
		Then the sequence $ (V_{N}(B ^{H}), N\geq 1) $ converges to $1$ as $N\to \infty$ in $ L^{p} (\Omega)$ for any $p\geq 1$  and almost surely.  Moreover, for $N$ large enough, and for every $p\geq 1$, 
		\begin{equation}
			\label{2f-1}
			\mathbf{E} \left| V _{N} (B ^{H})-1\right| ^{p} \leq C\left(  \frac{1}{N}\right) ^{\frac{p}{2}}.
		\end{equation}
	\end{prop}
	{\bf Proof: }The  bound (\ref{2f-1}) (and consequently the convergence in $ L ^{p}(\Omega)$) follows from the proof of Proposition 1 in \cite{GT}. Altough the almost sure convergence is not stated in \cite{GT}, it follows easily by using the Borel-Cantelli lemma. Indeed, for $\gamma \in (0, \frac{1}{2})$ and $p\geq 1$, we have by (\ref{2f-1}),
	\begin{eqnarray*}
		&&	\sum_{N \geq 1} P \left( \vert V_{N} (B ^{H}) -1\vert \geq N ^{-\gamma} \right) \leq \sum_{N \geq 1} N ^{p\gamma} 	\mathbf{E} \left| V_{N} (B ^{H}) -1\right| ^{p} \\
		&&\leq C \sum_{N \geq 1} N ^{p(\gamma-\frac{1}{2})}, 
	\end{eqnarray*} 
	and the last series converges for $p$ sufficiently large. \qed 
	
	We consider a perturbed fractional Brownian motion $ (X_{t}, t\geq 0)$, i.e. a centered Gaussian process  given by 
	\begin{equation}
		\label{pfBm}
		X_{t}= C_{0}B ^{H}_{t}+ Y_{t}, \hskip0.4cm t\geq 0,
	\end{equation}
	where $C_{0}>0$ and  $ (Y_{t}, t\geq 0)$ is a centered Gaussian process, self-similar of index $H\in (0,1)$,  satisfying 
	\begin{equation}
		\label{cy}
		\mathbf{E} \left| Y_{t}- Y_{s} \right| ^{2} \leq Cs^{2H-2},
	\end{equation}
	for every $0< s< t$. We will deduce some similar results to those in Lemma \ref{ll11} and Proposition \ref{pp2} for the process $X$. Let us start with the following lemma. Recall that $t_{i}, i=0,1,...,N$ are given by (\ref{ti}).  In what follows, we denote by $C$ a strictly positive constant that may change from one place to another.

	\begin{lemma}\label{ll2}
		Let $ H\in \left(0, 1\right)$ and let $X$ be given by (\ref{pfBm}).   Then  for every $i\geq 1$,

		\begin{equation}\label{2f-2}
			\mathbf{E} \left[ N ^{2H-1} \left( X _{t_{i+1}} - X  _{t_{i}} \right)^{2} -\frac{C_{0}^{2}}{N} \right] ^{2}\leq C \frac{1}{N^{2}},
		\end{equation}
		and	for any $p\geq 1$,
		\begin{equation}\label{8f-1}
			\mathbf{E} \left| X_{t_{i+1}} -X_{t_{i}} \right| ^{p} \leq C\frac{1}{N ^{pH}}.
		\end{equation}
		
		\noindent For every $i,j\geq 1$ with $i\not=j$, 
		
		\begin{eqnarray}
			&&	\mathbf{E}\left[ \left( N ^{2H-1} \left( X _{t_{i+1}} - X  _{t_{i}} \right)^{2} -\frac{C_{0}^{2}}{N} \right)\left( N ^{2H-1} \left( X _{t_{j+1}} - X  _{t_{j}} \right)^{2} -\frac{C_{0}^{2}}{N} \right)\right]\nonumber \\
			&&\leq C\frac{1}{N^{2}} (\rho_{H}^{2}(i-j)+ i ^{H-1}+j^{H-1}), \label{14s-1}
		\end{eqnarray}
		with
		\begin{equation}
			\label{ro}
			\rho_{H}(v)= \frac{1}{2}(\vert v+1\vert ^{2H}+\vert v-1\vert ^{2H}-2\vert v\vert ^{2H}), \hskip0.4cm v\in \mathbb{Z}.
		\end{equation}
		
	\end{lemma}
	{\bf Proof: } We have, by (\ref{cy}) and the self-similarity of $Y$, for every $i=1,...,N$,
	\begin{equation}\label{23m-2}
		\mathbf{E} \left| Y_{t_{i+1}} - Y_{t_{i}} \right|  ^{2}= N ^{-2H} 	\mathbf{E} \left| Y_{i+1}- Y_{i} \right| ^{2} \leq C N ^{-2H} i^{2H-2}. 
	\end{equation}
	To prove the bound (\ref{2f-2}), we write
	\begin{eqnarray}
		&&	\mathbf{E} \left[ N ^{2H-1} \left( X _{t_{i+1}} - X  _{t_{i}} \right)^{2} -\frac{C_{0} ^{2}}{N} \right] ^{2}\nonumber \\
		&=& 	\mathbf{E} \left[ C_{0}^{2}N ^{2H-1} \left( B ^{H} _{t_{i+1}} - B ^{H} _{t_{i}} \right)^{2} -\frac{C_{0} ^{2}}{N}  +2C_{0} N ^{2H-1} \left( B ^{H} _{t_{i+1}} - B ^{H} _{t_{i}} \right) \left( Y _{t_{i+1}} - Y_{t_{i}} \right)\right.\nonumber\\
		&&\left. +N ^{2H-1} \left( Y _{t_{i+1}} - Y_{t_{i}} \right)^{2} \right] ^{2}\nonumber\\
		&\leq & C\left[ 	\mathbf{E} \left[ N ^{2H-1} \left( B ^{H} _{t_{i+1}} - B ^{H} _{t_{i}} \right)^{2} -\frac{1}{N} \right] ^{2}+ N ^{4H-2} 	\mathbf{E} \left( B ^{H} _{t_{i+1}} - B ^{H} _{t_{i}} \right)^{2}  \left( Y_{t_{i+1}} - Y _{t_{i}} \right)^{2}\right.\nonumber\\
		&& \left.+ N ^{4H-2} 	\mathbf{E}\left( Y_{t_{i+1}} - Y _{t_{i}} \right)^{4} \right]\nonumber\\
		&\leq &  C\left[ N ^{-2} + N ^{4H-2} 	\mathbf{E} \left( B ^{H} _{t_{i+1}} - B ^{H} _{t_{i}} \right)^{2}  \left( Y_{t_{i+1}} - Y _{t_{i}} \right)^{2} +  N ^{4H-2}	\mathbf{E}\left( Y_{t_{i+1}} - Y _{t_{i}} \right)^{4} \right],\label{2f-4}
	\end{eqnarray}
	where we used Lemma \ref{ll11}.	 Now,  via Cauchy-Schwarz, the property of the Gaussian moments and (\ref{23m-2}), we have for $i\geq 1$,
	\begin{eqnarray}
		&&N ^{4H-2}		\mathbf{E} \left( B ^{H} _{t_{i+1}} - B ^{H} _{t_{i}} \right)^{2}  \left( Y_{t_{i+1}} - Y _{t_{i}} \right)^{2} \nonumber\\
		&& \leq N ^{4H-2} \left( 	\mathbf{E} \left( B ^{H}_{t_{i+1}}-B^{H} _{t_{i}}\right) ^{4} \right) ^{\frac{1}{2}} \left( 	\mathbf{E} \left( Y _{t_{i+1}}-Y _{t_{i}}\right) ^{4} \right) ^{\frac{1}{2}}\nonumber\\
		&&=C N ^{4H-2} 	\mathbf{E} \left( B ^{H}_{t_{i+1}}-B^{H} _{t_{i}}\right) ^{2}	\mathbf{E} \left( Y_{t_{i+1}}-Y _{t_{i}}\right) ^{2}\leq  C N ^{-2} i ^{2H-2}\leq C N ^{-2},\label{2f-5}
	\end{eqnarray}
	and
	\begin{eqnarray}
		&&	N ^{4H-2}	\mathbf{E}\left( Y_{t_{i+1}} - Y _{t_{i}} \right)^{4} \leq CN ^{4H-2} \left( 	\mathbf{E}\left( Y_{t_{i+1}} - Y _{t_{i}} \right)^{2} \right) ^{2} \nonumber\\
		&&\leq C N ^{-2} i ^{4H-4}\leq C N ^{-2}.\label{2f-6}
	\end{eqnarray}
	By plugging (\ref{2f-5}) and (\ref{2f-6}) into (\ref{2f-4}), we get the bound (\ref{2f-2}). Finally, for $N\geq 1$ and $i=1,...,N$,
	\begin{eqnarray*}
		\mathbf{E}  \left| X_{t_{i+1}} -X_{t_{i}} \right| ^{p} &=&\sum_{k=0} ^{p} C_{p} ^{k} C_{0} ^{p-k} 	\mathbf{E}\left(   \left| B ^{H}_{t_{i+1}} -B ^{H}_{t_{i}} \right| ^{p-k}  \left| Y_{t_{i+1}} -Y_{t_{i}} \right|  ^{k} \right)\\
		&\leq & \sum_{k=0} ^{p} C_{p} ^{k} C_{0} ^{p-k} \left( 	\mathbf{E}  \left| B ^{H}_{t_{i+1}} -B ^{H}_{t_{i}} \right| ^{2p-2k}\right) ^{\frac{1}{2}} \left( 	\mathbf{E}  \left| Y_{t_{i+1}} -Y_{t_{i}} \right|  ^{2k}\right) ^{\frac{1}{2}}\\
		&\leq & C\frac{1}{N ^{pH}},
	\end{eqnarray*}
	and thus (\ref{8f-1}) is obtained.  
	
	Let us finally show (\ref{14s-1}). We write
	\begin{eqnarray}
		&&	\mathbf{E}\left[ \left( N ^{2H-1} \left( X _{t_{i+1}} - X  _{t_{i}} \right)^{2} -\frac{C_{0}^{2}}{N} \right)\left( N ^{2H-1} \left( X _{t_{j+1}} - X  _{t_{j}} \right)^{2} -\frac{C_{0}^{2}}{N} \right)\right]\nonumber\\
		&=&	\mathbf{E}\left[ \left( N ^{2H-1}C_{0}^{2} \left( B^{H} _{t_{i+1}} - B^{H}  _{t_{i}} \right)^{2} -\frac{C_{0}^{2}}{N} \right.\right.\nonumber\\
		&&\left. \left.+2C_{0}N ^{2H-1}\left( B^{H} _{t_{i+1}} - B^{H}  _{t_{i}} \right)\left( Y _{t_{i+1}} - Y  _{t_{i}} \right)+ N ^{2H-1} \left( Y_{t_{i+1}} -Y_{t_{i}}\right) ^{2}\right)\right.\nonumber\\
		&&\left. \left( N ^{2H-1}C_{0}^{2} \left( B^{H} _{t_{j+1}} - B^{H}  _{t_{j}} \right)^{2} -\frac{C_{0}^{2}}{N} \right.\right.\nonumber\\
		&&\left. \left. +2C_{0}N ^{2H-1}\left( B^{H} _{t_{j+1}} - B^{H}  _{t_{j}} \right)\left( Y _{t_{j+1}} - Y  _{t_{j}} \right)+ N ^{2H-1} \left( Y_{t_{j+1}} -Y_{t_{j}}\right) ^{2}\right)\right]. \label{14s-2}
	\end{eqnarray}
	By Lemma 3 in \cite{GT},
	\begin{eqnarray*}
		&&	\mathbf{E}\left[ \left( N ^{2H-1}C_{0}^{2} \left( B^{H} _{t_{i+1}} - B^{H}  _{t_{i}} \right)^{2} -\frac{C_{0}^{2}}{N} \right) \left( N ^{2H-1}C_{0}^{2} \left( B^{H} _{t_{j+1}} - B^{H}  _{t_{j}} \right)^{2} -\frac{C_{0}^{2}}{N} \right)\right]\\
		&&\leq C \rho_{H}^{2} (i-j)\frac{1}{N ^{2}},
	\end{eqnarray*}
	with $\rho_{H}$ given by (\ref{ro}). Also,
	\begin{eqnarray*}
		&&\mathbf{E}\left[ \left( N ^{2H-1}C_{0}^{2} \left( B^{H} _{t_{i+1}} - B^{H}  _{t_{i}} \right)^{2} -\frac{C_{0}^{2}}{N} \right) N ^{2H-1}\left( B^{H} _{t_{j+1}} - B^{H}  _{t_{j}} \right)\left( Y _{t_{j+1}} - Y  _{t_{j}} \right)\right] \\
		&\leq & C N ^{2H-1} \left( \mathbf{E}  \left( N ^{2H-1}C_{0}^{2} \left( B^{H} _{t_{i+1}} - B^{H}  _{t_{i}} \right)^{2} -\frac{C_{0}^{2}}{N} \right)^{2}\right) ^{\frac{1}{2}} \\
		&&\times \left( \mathbf{E}\left( B^{H} _{t_{j+1}} - B^{H}  _{t_{j}} \right)^{4} \right)^{\frac{1}{4}} \left( \mathbf{E}\left(Y _{t_{j+1}} - Y _{t_{j}} \right)^{4} \right)^{\frac{1}{4}}\\
		&\leq & C\frac{ j ^{H-1}}{N^{2}},
	\end{eqnarray*}
	where we used Lemma \ref{ll11}  and (\ref{23m-2}) for the last line. Similarly,
	\begin{equation*}
		\mathbf{E}\left[ \left( N ^{2H-1}C_{0}^{2} \left( B^{H} _{t_{j+1}} - B^{H}  _{t_{j}} \right)^{2} -\frac{C_{0}^{2}}{N} \right) N ^{2H-1}\left( B^{H} _{t_{i+1}} - B^{H}  _{t_{i}} \right)\left( Y _{t_{i+1}} - Y  _{t_{i}} \right)\right] \leq C\frac{ i^{H-1}}{N^{2}}
	\end{equation*}
	and 
	\begin{eqnarray*}
		&&	\mathbf{E}\left[ \left( N ^{2H-1}C_{0}^{2} \left( B^{H} _{t_{j+1}} - B^{H}  _{t_{j}} \right)^{2} -\frac{C_{0}^{2}}{N} \right) N ^{2H-1}\left( Y_{t_{i+1}} - Y _{t_{i}} \right)^{2} \right] \\
		&&\leq  N ^{2H-1}\left( \mathbf{E} \left( N ^{2H-1}C_{0}^{2} \left( B^{H} _{t_{j+1}} - B^{H}  _{t_{j}} \right)^{2} -\frac{C_{0}^{2}}{N} \right) ^{2} \right) ^{\frac{1}{2}} \left( \mathbf{E} \left( Y_{t_{i+1}} - Y _{t_{i}} \right)^{4} \right) ^{\frac{1}{2}}\\
		&&\leq C N ^{2H-1} \frac{1}{N} N ^{-2H} i ^{2H-2} \leq C \frac{ i^{2H-2}}{N ^{2}}
	\end{eqnarray*}
	and
	\begin{eqnarray*}
		&&	N ^{4H-2} \mathbf{E}  \left( B^{H} _{t_{i+1}} - B^{H}  _{t_{i}}\right) \left( Y _{t_{i+1}} - Y _{t_{i}}\right) \left( Y _{t_{j+1}} -Y  _{t_{j}}\right) \\
		&&\leq N ^{4H-2}  \left( \mathbf{E} \left( B^{H} _{t_{i+1}} - B^{H}  _{t_{i}}\right) ^{4}  \right) ^{\frac{1}{4}} \left( \mathbf{E} \left(Y_{t_{i+1}} - Y  _{t_{i}}\right) ^{4}  \right) ^{\frac{1}{4}}\left( \mathbf{E} \left(Y_{t_{i+1}} - Y  _{t_{i}}\right) ^{4}  \right) ^{\frac{1}{2}}\\
		&&\leq C N^{4H-2} N ^{-H} i^{H-1} N ^{-H} j^{2H-2} N ^{-2H}\leq C \frac{ i ^{H-1}}{N ^{2}}.
	\end{eqnarray*}
	The other summands appearing when we develop the right-hand side of (\ref{14s-2}) can be treated in a similar manner. 
	\qed

	We can immediately deduce the behavior of the quadratic variation of the perturbed fBm.

	\begin{prop}\label{pp22}
		Let $( B ^{H}_{t}, t\geq 0)$ be a fBm with $H\in \left(0, \frac{3}{4}\right)$  and let $ (X_{t}, t\geq 0) $ be a perturbed fBm defined by (\ref{pfBm}). For $N \geq 1$, let
		\begin{equation*}
			V_{N}(X) = N ^{2H-1} \sum_{i=0} ^{N-1} \left( X_{\frac{i+1}{N}} -X _{\frac{i}{N}} \right)^{2}.
		\end{equation*}
		Then the sequence $(V_{N}(X), N\geq 1)$ converges to $ C_{0} ^{2} $ in $ L ^{p} (\Omega)$ for every $p\geq 1$ and almost surely, with $ C_{0} $ from (\ref{pfBm}). Moreover, for $N$ large enough, 
		\begin{equation}\label{2f-10}
			\mathbf{E} \left| V_{N}(X) - C_{0} ^{2} \right| ^{2} \leq C \begin{cases}
				N ^{-1} \mbox{ if } H \in (0, \frac{1}{2}),\\
				\log(N) N ^{-1} \mbox{ if } H=\frac{1}{2},\\
				N ^{2H-2}, \mbox{ if } H \in (\frac{1}{2}, \frac{3}{4}).
			\end{cases}
		\end{equation}
	\end{prop}
	{\bf Proof: } We can write, for every $N\geq 1$, 
	\begin{eqnarray*}
		V_{N} (X) &=& C_{0}^{2} V_{N} (B ^{H})+ 2C_{0} N ^{2H-1}  \sum_{i=0} ^{N-1} \left( B ^{H}_{t_{i+1}} -B ^{H} _{t_{i}} \right)\left( Y_{t_{i+1}} -Y _{t_{i}} \right)\\
		&&+ N ^{2H-1} \sum_{i=0} ^{N-1} \left( Y_{t_{i+1}} -Y _{t_{i}} \right)^{2}.
	\end{eqnarray*}
	Thus,
	\begin{eqnarray}
		&&	\mathbf{E} \left| V_{N}(X) - C_{0} ^{2} \right| ^{2}\nonumber  \\
		&&\leq C \left[ 	\mathbf{E} \left| V_{N}(B ^{H}) - 1\right| ^{2} + N ^{4H-2} 	\mathbf{E} \left| \sum_{i=0} ^{N-1} \left( B^{H}_{ t_{i+1}}-B ^{H} _{t_{i}}\right) \left( Y_{ t_{i+1}}-Y _{t_{i}}\right) \right| ^{2}\right.\nonumber\\
		&&\left.  + N ^{4H-2} 	\mathbf{E} \left| \sum_{i=0} ^{N-1} \left( Y_{ t_{i+1}}-Y _{t_{i}}\right) ^{2} \right| ^{2} \right]. \label{2f-7}
	\end{eqnarray}

	By Proposition \ref{pp2}, $	\mathbf{E} \left| V_{N}(B ^{H}) - 1\right| ^{2} $ converges to 0 as $N\to \infty$ and the bound (\ref{2f-1}) holds true. Next,  by (\ref{23m-2}), 
	
	\begin{eqnarray}
		&&N ^{4H-2} 	\mathbf{E} \left| \sum_{i=0} ^{N-1} \left( Y_{ t_{i+1}}-Y _{t_{i}}\right) ^{2} \right| ^{2}\nonumber\\
		&&\leq N ^{4H-1} \sum_{i=0} ^{N-1} 	\mathbf{E}  \left( Y_{ t_{i+1}}-Y _{t_{i}}\right) ^{4}\leq C N ^{4H-1} \sum_{i=0} ^{N-1} \left( 	\mathbf{E}  \left( Y_{ t_{i+1}}-Y _{t_{i}}\right) ^{2}\right) ^{2}\nonumber\\
		&&\leq C N ^{4H-1}\left(1+ \sum_{i=1} ^{N-1} \left( 	\mathbf{E}  \left( Y_{ t_{i+1}}-Y _{t_{i}}\right) ^{2}\right) ^{2}\right)\nonumber\\
		&&\leq C \frac{1}{N}\left( 1+ \sum_{i=1} ^{N-1} i ^{4H-4} \right)\leq C \frac{1}{N},\label{2f-8}
	\end{eqnarray}
	because $\sum_{i\geq 1} i ^{4H-4}$ converges for $H<\frac{3}{4}$. Also,

	\begin{eqnarray}
		&&	N ^{4H-2} 	\mathbf{E} \left| \sum_{i=0} ^{N-1} \left( B^{H}_{ t_{i+1}}-B ^{H} _{t_{i}}\right) \left( Y_{ t_{i+1}}-Y _{t_{i}}\right) \right| ^{2}\nonumber\\
		&&\leq N ^{4H-1} \sum_{i=0} ^{N-1} 	\mathbf{E}  \left( B^{H}_{ t_{i+1}}-B ^{H} _{t_{i}}\right) ^{2} \left( Y_{ t_{i+1}}-Y _{t_{i}}\right) ^{2}\nonumber \\
		&&\leq C N ^{4H-1}\sum_{i=0} ^{N-1} 	\mathbf{E}  \left( B^{H}_{ t_{i+1}}-B ^{H} _{t_{i}}\right) ^{2}	\mathbf{E} \left( Y_{ t_{i+1}}-Y _{t_{i}}\right) ^{2}\nonumber\\
		&&\leq C \frac{1}{N}\left(1+ \sum_{i=1} ^{N-1} i ^{2H-2}\right)\nonumber\\
		&&\leq C \begin{cases}
			N ^{-1} \mbox{ if } H \in (0, \frac{1}{2}),\\
			\log(N) N ^{-1} \mbox{ if } H=\frac{1}{2},\\
			N ^{2H-2}, \mbox{ if } H \in (\frac{1}{2}, \frac{3}{4}).
		\end{cases}\label{2f-9}
	\end{eqnarray}
	We combine (\ref{2f-8}) and (\ref{2f-9}) with (\ref{2f-7}) and we get the inequality (\ref{2f-10}).  This implies the convergence in $ L^{2}(\Omega)$ and thus the convergence in $ L^{p}(\Omega)$ for every $p\geq 1$, due to the hypercontractivity property  (see e.g. \cite{NP-book}), since $ V_{N}(X) $ is an element in a finite sum of Wiener chaoses  and then, for every $p\geq 1$, 
	\begin{equation}\label{11f-1}
		\mathbf{E} \left| V_{N}(X) - C_{0} ^{2} \right| ^{p} \leq C_{p} \left(	\mathbf{E} \left| V_{N}(X) - C_{0} ^{2} \right| ^{2} \right) ^{\frac{p}{2}}	\end{equation}
	
	The almost sure convergence is obtained via Borel-Cantelli, as in the proof of Proposition \ref{pp2}.  \qed

	\section{Quadratic and higher order variations of the solution to the nonlinear stochastic heat equation}
	In this paragraph, we study the asymptotic behavior of the renormalized quadratic variation and of another $p$-variation (with some particular $p>2$) of the solution to the stochastic heat equation (\ref{2}). These variations are considered with respect to the time variable. We start by obtaining the case of the  heat equation with additive noise and then we extend our results to the nonlinear noise case. The limit theorems obtained below will be used in the next section for statistical inference.

	\subsection{The (renormalized) quadratic variation of the solution to the fractional heat equation with additive space-time white noise}
	
	In a first step, we deduce the limit behavior of the quadratic variation of the solution in the additive noise case, i.e. $\sigma (x)= 1$ for every $x \in \mathbb{R}$ in (\ref{2}). We denote by $u_{0}$ the solution to this equation, i.e. for every $t\geq 0, x\in \mathbb{R}$,
	\begin{equation}
		\label{u0}
		u_{0}(t,x)=\int_{0} ^{t} \int _{\mathbb{R}}G_{\alpha} (t-s, x-y) W(ds, dy),
	\end{equation}
	where $ G_{\alpha}$ is the Green kernel given by (\ref{ga}). Then $ (u_{0}(t,x), t\geq 0, x\in \mathbb{R})$ is a centered Gaussian  random field and viewed as a stochastic process with respect to its temporal variable, $u_{0}$ is related to the fractional Brownian motion in the following way  (see e.g. \cite{T2}): 
	\begin{equation}
		\label{deco}
		\left( u_{0}(t,x), t\geq 0\right)= \left( C_{0, \alpha} U_{t}+ Y_{t}, t\geq 0\right),
	\end{equation} 
	where $C_{0, \alpha}>0$ (its explicit expression can be found in e.g. \cite{T2}),  $ (U_{t}, t\geq 0)$ is a fractional Brownian motion with Hurst index $H=\frac{\alpha-1}{2\alpha}$ and $(Y_{t}, t\geq 0)$ is a $H$-self-similar Gaussian process satisfying (\ref{cy}), i.e. for $i\geq 1$ ,
	\begin{equation}\label{12f-1}
		\mathbf{E} \left| Y_{i+1}- Y_{i} \right| ^{2} \leq  C i ^{-\frac{1+\alpha }{\alpha}}.
	\end{equation}
	Notice that $H\in \left(0, \frac{1}{4}\right) $ if $\alpha \in (1,2)$ and $H=\frac{1}{4}$ if $\alpha =2$. The bound (\ref{12f-1}) can be deduced from Lemma 3.5 in \cite{BEV} or Section 5.2.2 in \cite{T2}. In other words, the process $u_{0}$ is, with respect to its time variable, a perturbed fractional Brownian motion in the sense of the definition (\ref{pfBm}). The same happens for $u_{0}$ with respect to its space variable, but the perturbation $Y$ has different properties (see e.g. \cite{GT}).
	
	By (\ref{12f-1}) and the scaling property, for $N\geq 1$ and $i=1,..., N-1$, 
	\begin{equation}\label{23i-1}
		\mathbf{E} \left| Y_{t_{i+1}}- Y_{t_{i}} \right| ^{2}=\left( \frac{1}{N}\right) ^{\frac{ \alpha-1}{\alpha}}\mathbf{E} \left| Y_{i+1}- Y_{i} \right| ^{2} \leq C \left( \frac{1}{N}\right) ^{\frac{ \alpha-1}{\alpha}}i ^{-\frac{2}{\alpha}}.
	\end{equation}
	Let us set, for $N\geq 1$ and for $x\in \mathbb{R}$,
	\begin{equation}\label{vnxu0}
		V_{N,x}(u_{0})=N ^{-\frac{1}{\alpha}}\sum_{i=0} ^{N-1} \left( u_{0} (t_{i+1},x)- u_{0} (t_{i}, x) \right) ^{2},
	\end{equation}
	with $t_{i}, i=0,1,...,N$ given by (\ref{ti}). The limit behavior of the above sequence is an immediate consequence of Proposition  \ref{pp2}.
	
	\begin{prop}\label{pp3}
		Let $x\in \mathbb{R}$ be fixed and consider the sequence $ (V_{N,x} (u_{0}), N\geq 1)$ given by (\ref{vnxu0}). Then 
		\begin{equation*}
			V_{N,x}(u_{0})\to _{N \to \infty} C_{0, \alpha} ^{2} \mbox{almost surely and in } L ^{2}(\Omega).
		\end{equation*}
		with $C_{0, \alpha}$ from (\ref{deco}).  For $N$ large enough, one has
		\begin{equation*}
			\mathbf{E} \left| V_{N,x}(u_{0})- C_{0, \alpha}^{2} \right| ^{2} \leq C \frac{1}{N}.
		\end{equation*}
	\end{prop}
	{\bf Proof: } Via (\ref{deco}) and the inequality (\ref{23i-1}), it suffices to apply Proposition \ref{pp22} with $H=\frac{\alpha-1}{2\alpha}$, which actually belongs to the interval $\left(0, \frac{1}{2}\right). $ \qed

	\subsection{The nonlinear case}
	The purpose is to obtain a similar result to Proposition 3 in the case of the nonlinear heat equation. The idea is to compare the temporal  increments of the process (\ref{mild1}) with those of the solution $u_{0}$  to the linear heat equation. 
	For $\delta>0$, we denote by $\Delta_{1}u_{0}(t, \delta)$ the temporal increment of lenght $\delta$ of the random field $u_{0}$, i.e. 
	\begin{eqnarray*}
		\Delta_{1}u_{0}(t, \delta)&=& u_{0} (t+\delta, x)- u_{0}(t,x)\\
		&=& \int_{0} ^{t+\delta} \int_{\mathbb{R}} \left( G_{\alpha}(t+\delta -a, x-y)- G_{\alpha}(t-a, x-y)\right)W(da, dy).
	\end{eqnarray*}
	Recall that $ G_{\alpha}(t,x)=0$ if $t\leq 0$. 
	
	Let $\widetilde{W}$ be an independent copy of the space-time white noise $W$ and let us consider the modified temporal increment of $u_{0}$ given by 
	\begin{eqnarray*}
		\widetilde{\Delta}_{1} u_{0} (t,\delta) &=& \int_{0} ^{t(\delta)}\int_{\mathbb{R}} \left( G_{\alpha}(t+\delta -a, x-y)- G_{\alpha}(t-a, x-y)\right)\widetilde{W}(da, dy)\\
		&&+ \int_{ t(\delta)} ^{t+\delta} \int_{\mathbb{R}} \left( G_{\alpha}(t+\delta -a, x-y)- G_{\alpha}(t-a, x-y)\right)W(da, dy),
	\end{eqnarray*} 
	where we set
	\begin{equation}
		\label{td}
		t(\delta)= t-\delta  ^{\beta} \mbox{ with } \beta=\frac{2\alpha}{2\alpha+1},
	\end{equation}
	if $ \delta ^{\beta} <t$. 	For every $t\geq 0, x\in \mathbb{R}$ and $\delta >0$, we have:
	
	\begin{enumerate}
		\item The random variable $\widetilde{\Delta}_{1} u_{0} (t,\delta) $ has the same distribution as the standard increment $	\Delta_{1}u_{0}(t, \delta)$. In particular, it is a centered Gaussian random variable. 
		
		\item The random variable $\widetilde{\Delta}_{1} u_{0} (t,\delta) $ is independent of the sigma-algebra $\mathcal{G}= \sigma \{ u(s,x), 0\leq s\leq t(\delta), x\in \mathbb{R}\}$. 
	\end{enumerate}
	
	The below result plays an important role in the sequel. It allows to compare the increments of the solution to the nonlinear heat equation with those of the solution to the stochastic heat equation with additive space-time white noise. 
	
	\begin{prop}\label{pp4}
		For every $x\in \mathbb{R}$ and $\delta >0$,
		\begin{equation}\label{1a-1}
			\mathbf{E} \left| 	\left( u(t+\delta, x)- u(t,x)\right) - \sigma (u( t(\delta), x))( \widetilde{\Delta} _{1} u_{0} ) (t, \delta) \right| ^{2}\leq C\delta ^{ \frac{ 4(\alpha -1)}{2\alpha +1}}.
		\end{equation}
	\end{prop}
	{\bf Proof: } The proof is postponed to the Appendix. \qed 
	
	If $\alpha =2$, the estimate in the right-hand side of (\ref{1a-1}) reads $C \delta ^{\frac{4}{5}}$, which corresponds with the bound obtained in Lemma 4.5 in \cite{PoTr}.
	Let us now define
	\begin{equation}
		\label{vnxu}
		V_{N,x}(u)= N ^{-\frac{1}{\alpha}} \sum_{i=0} ^{N-1} \left( u(t_{i+1}, x)-u(t_{i}, x)\right) ^{2}, \hskip0.4cm N\geq 1,
	\end{equation}
	where $u$ is the mild solution (\ref{mild1}) and $t_{i}, i=0,1...,N$ are given by (\ref{ti}). The sequence (\ref{vnxu}) constitutes the renormalized temporal  quadratic variation of the process  $u$ given by (\ref{mild1}).

	We have the following result.
	
	\begin{theorem}\label{tt1}
		Let $x\in \mathbb{R}$ be fixed and consider the sequence $(V_{N,x}(u), N\geq 1)$ given by (\ref{vnxu}). Then, for $N$ sufficiently large,
		\begin{equation}
			\label{3f-5}	\mathbf{E} \left| 	V_{N,x}(u)-C_{0, \alpha} ^{2} \int_{0}^{1} \sigma ^{2} (u(s,x))ds\right| \leq C N^{\frac{ (\alpha -1)(-2\alpha +1)}{2\alpha (2\alpha +1)}},
		\end{equation}
		with $C_{0, \alpha}$ from (\ref{deco}). In particular,
		\begin{equation*}
			V_{N,x}(u)\to _{N\to \infty} C_{0, \alpha} ^{2} \int_{0}^{1} \sigma ^{2} (u(s,x))ds \mbox{ in } L ^{1} (\Omega). 
		\end{equation*} 
		
	\end{theorem}
	{\bf Proof: } We write
	\begin{eqnarray*}
		&&	V_{N,x}(u)-C_{0, \alpha}^{2}  \int_{0}^{1} \sigma ^{2} (u(s,x))ds\\
		&=& N ^{-\frac{1}{\alpha}}\left(\sum_{i=0} ^{N-1} \left( \left( u(t_{i+1}, x)-u(t_{i}, x)\right) ^{2}-\sigma ^{2} (u(t_{i}(\delta), x))(\widetilde{\Delta} _{1} u_{0} (t_{i}, \delta))^{2}\right)\right) \\
		&&+\sum_{i=0} ^{N-1} \sigma ^{2} (u(t_{i}(\delta), x)) \left( N ^{-\frac{1}{\alpha}}(\widetilde{\Delta} _{1} u_{0} (t_{i}, \delta))^{2}- C_{0, \alpha}^{2} \frac{1}{N}\right) \\
		&&+C_{0, \alpha } ^{2} \left[ \frac{1}{N} \sum_{i=0} ^{N-1} \sigma ^{2} (u(t_{i}(\delta), x))- \int_{0}^{1} \sigma ^{2} (u(s,x))ds\right]\\
		&:=& L_{1,N}+ L_{2,N}+ L_{3,N}.
	\end{eqnarray*}
	Let us estimate the three terms $L_{1,N}, L_{2,N}$ and $ L_{3,N}$ when $N$ is large enough. To deal with $ L_{1,N}$, we use Proposition \ref{pp4} with $t=t_{i}$ and $\delta =\frac{1}{N}$. We get 
	\begin{eqnarray*}
		\mathbf{E} \vert L_{1, N} \vert &\leq&  N ^{-\frac{1}{\alpha}} \sum_{i=0} ^{N-1}
		\left( 	\mathbf{E}	\left|  \left( u(t_{i+1}, x) - u(t_{i}, x) \right) -\sigma  (u(t_{i}(\delta),x)(\widetilde{\Delta} _{1} u_{0} )(t_{i}, \delta)  \right| ^{2}  \right) ^{\frac{1}{2}}\\
		&&\times \left( 	\mathbf{E}	\left|  \left( u(t_{i+1}, x) - u(t_{i}, x) \right) +\sigma  (u(t_{i}(\delta),x)(\widetilde{\Delta} _{1} u_{0} )(t_{i}, \delta)  \right|^{2} \right) ^{\frac{1}{2}}\\
		&\leq & C   N ^{-\frac{1}{\alpha}}  \left( \frac{1}{N}\right) ^{\frac{2(\alpha -1)}{2\alpha +1}}\\
		&& \sum_{i=0} ^{N-1}
		\left(	\mathbf{E}	\left|  \left( u(t_{i+1}, x) - u(t_{i}, x) \right) +\sigma  (u(t_{i}(\delta),x))(\widetilde{\Delta} _{1} u_{0} )(t_{i}, \delta)  \right|  ^{2} \right) ^{\frac{1}{2}}.
	\end{eqnarray*}
	We notice that, by (\ref{2i-2}) and the comment following (\ref{td}),
	\begin{equation}
		\label{17s-1}
		\mathbf{E}\vert  (\widetilde{\Delta} _{1} u_{0} )(t_{i}, \delta) \vert ^{2}\leq C N ^{-\frac{\alpha -1}{2\alpha}}.
	\end{equation}
	Now, via (\ref{2i-2}), (\ref{2i-3}) and (\ref{17s-1}), we get the bound  
	\begin{eqnarray*}
		&&	\mathbf{E}	\left|  \left( u(t_{i+1}, x) - u(t_{i}, x) \right) +\sigma  (u(t_{i}(\delta),x)(\widetilde{\Delta} _{1} u_{0} )(t_{i}, \delta)  \right|^{2} \\
		&&\leq C \left( \mathbf{E}  \left( u(t_{i+1}, x) - u(t_{i}, x) \right) ^{2} + \mathbf{E} \sigma  (u(t_{i}(\delta),x)^{2} \mathbf{E} (\widetilde{\Delta} _{1} u_{0} )(t_{i}, \delta) ^{2} \right) \\
		&&\leq  C \left( \mathbf{E}  \left( u(t_{i+1}, x) - u(t_{i}, x) \right) ^{2} +  \mathbf{E} (\widetilde{\Delta} _{1} u_{0} )(t_{i}, \delta) ^{2} \right) \\
		&&\leq C N ^{ -\frac{\alpha -1}{\alpha}}.
	\end{eqnarray*}
	Hence
	\begin{equation}\label{3f-11}
		\mathbf{E} \vert L_{1, N} \vert \leq C   N ^{-\frac{1}{\alpha}}  \left( \frac{1}{N}\right) ^{\frac{2(\alpha -1)}{2\alpha +1}}N ^{ -\frac{\alpha -1}{2\alpha}+1}=C  N^{\frac{ (\alpha -1)(-2\alpha +1)}{2\alpha (2\alpha +1)}}\to _{N \to \infty} 0. 
	\end{equation}
	We estimate  $ L_{2,N}$ as follows
	\begin{eqnarray*}
		\mathbf{E} L_{2,N}^{2} &=& 	\mathbf{E} \sum_{i,j=0} ^{N-1} \sigma^{2}   (u(t_{i}(\delta),x)\sigma  ^{2} (u(t_{j}(\delta),x)\\
		&&\left( N ^{-\frac{1}{\alpha}}(\widetilde{\Delta} _{1} u_{0} (t_{i}, \delta))^{2}- C_{0, \alpha}^{2} \frac{1}{N}\right)\left( N ^{-\frac{1}{\alpha}}(\widetilde{\Delta} _{1} u_{0} (t_{j}, \delta))^{2}- C_{0, \alpha}^{2} \frac{1}{N}\right)\\
		&=& \sum_{i=0} ^{N-1} 	\mathbf{E}  \sigma^{4}   (u(t_{i}(\delta),x)	\mathbf{E}\left( N ^{-\frac{1}{\alpha}}(\widetilde{\Delta} _{1} u_{0} (t_{i}, \delta))^{2}- C_{0, \alpha}^{2} \frac{1}{N}\right)^{2} \\
		&&+ 2	\mathbf{E} \sum_{i,j=0;i>j} ^{N-1} \sigma^{2}   (u(t_{i}(\delta),x)\sigma  ^{2} (u(t_{j}(\delta),x)\\
		&&\left( N ^{-\frac{1}{\alpha}}(\widetilde{\Delta} _{1} u_{0} (t_{i}, \delta))^{2}- C_{0, \alpha}^{2} \frac{1}{N}\right)\left( N ^{-\frac{1}{\alpha}}(\widetilde{\Delta} _{1} u_{0} (t_{j}, \delta))^{2}- C_{0, \alpha}^{2} \frac{1}{N}\right)\\
		&:=& l_{2,N}^{(1)}+ l_{2,N}^{(2)},
	\end{eqnarray*}
	where we used the fact that $\widetilde{\Delta} _{1} u_{0} (t_{i}, \delta)$ is independent of $ \sigma \{ u(s, x), s\leq t_{i}, x\in \mathbb{R}\}$.  By the Lipschitz assumption on $\sigma$, (\ref{3i-5}) and the bound (\ref{2f-2}) in Lemma \ref{ll2},
	\begin{eqnarray*}
		l_{2,N}^{(1) }&\leq & C \sum_{i=0} ^{N-1}	\mathbf{E} \left( N ^{-\frac{1}{\alpha}}(\widetilde{\Delta} _{1} u_{0} (t_{i}, \delta))^{2}- C_{0, \alpha}^{2} \frac{1}{N}\right)^{2} \\
		&\leq & C N ^{-1}.
	\end{eqnarray*}
	In order to deal with the summande $ l_{2,N}^{(2)}$, we will consider the following situations:
	
	\begin{itemize}
		\item When $ i-j \leq [N ^{-\beta+1}]+1$.  The corresponding sum will be denoted by $ l_{2, N} ^{(2,1)}$.
		
		\item  When  $ i-j \geq [N ^{-\beta+1}]+2$. In this case we have
		$$t_{j}+ \delta \leq t_{i} (\delta).$$
		This implies that  $\widetilde{\Delta}_{1} u_{0} (t_{i}, \delta)$  and $\widetilde{\Delta}_{1} u_{0} (t_{j}, \delta)$  are independent of  $\sigma^{2}   (u(t_{i}(\delta),x)$ and $\sigma  ^{2} (u(t_{j}(\delta),x)$. The corresponding sum will de denoted by $ l_{2,N} ^{(2,2)}$. 
	\end{itemize}
	We denoted by $[x]$ the biggest integer less or equal than $x$. 	Let us now estimate the quantities $ l_{2,N} ^{(2,1)}$ and $ l_{2, N} ^{(2,2)}$. We can write, by Cauchy-Schwarz, the Lipschitz condition (\ref{2i-3}),  (\ref{3i-5})  and the hypercontractivity property (\ref{11f-1}),
	\begin{eqnarray*}
		l_{2, N} ^{(2,1)}&\leq & \sum_{i,j=0; i-j\leq [N ^{-\beta +1}]+1}^{N-1} \left( 	\mathbf{E}  \sigma^{8}   (u(t_{i}(\delta),x))\right) ^{\frac{1}{4}}  \left( 	\mathbf{E}  \sigma^{8}   (u(t_{j}(\delta),x))\right) ^{\frac{1}{4}} \\&&
		\times  \left( 	\mathbf{E} \left( N ^{-\frac{1}{\alpha}}(\widetilde{\Delta} _{1} u_{0} (t_{i}, \delta))^{2}- C_{0, \alpha}^{2} \frac{1}{N}\right)^{4} \right) ^{\frac{1}{4}} \\
		&&\times  \left( 	\mathbf{E} \left( N ^{-\frac{1}{\alpha}}(\widetilde{\Delta} _{1} u_{0} (t_{j}, \delta))^{2}- C_{0, \alpha}^{2} \frac{1}{N}\right)^{4} \right) ^{\frac{1}{4}}\\
		&\leq & C \sum_{i,j=0; i-j\leq [N ^{-\beta +1}]+1}^{N-1} \left( 	\mathbf{E} \left( N ^{-\frac{1}{\alpha}}(\widetilde{\Delta} _{1} u_{0} (t_{i}, \delta))^{2}- C_{0, \alpha}^{2} \frac{1}{N}\right)^{2} \right) ^{\frac{1}{2}} \\
		&&\times  \left( 	\mathbf{E} \left( N ^{-\frac{1}{\alpha}}(\widetilde{\Delta} _{1} u_{0} (t_{j}, \delta))^{2}- C_{0, \alpha}^{2} \frac{1}{N}\right)^{2} \right) ^{\frac{1}{2}}.
	\end{eqnarray*}
	By using the inequality (\ref{2f-2}) in Lemma \ref{ll2}, we obtain
	
	\begin{equation*}
		l_{2, N} ^{(2,1)}\leq C \frac{1}{ N ^{2}} \sum_{i,j=0; i-j\leq [N ^{-\beta +1}]+1}^{N-1} \leq C N ^{-\beta} = C N ^{-\frac{2\alpha}{2\alpha +1}}.
	\end{equation*}
	For the summand denoted by $l_{2,N} ^{(2,2)}$, we use the fact that  $\widetilde{\Delta}_{1} u_{0} (t_{i}, \delta)$ and $\widetilde{\Delta}_{1} u_{0} (t_{j}, \delta)$ are independent  of the other random variables appearing under the sum symbol, and we can write, via (\ref{2i-3}) and (\ref{3i-5}),
	
	\begin{eqnarray*}
		l_{2,N} ^{(2,2)} 
		&=& 2\sum_{i,j=0; i-j\geq [N ^{-\beta +1}]+2}^{N-1} 
		\mathbf{E} \left[ \sigma^{2}   (u(t_{i}(\delta),x))\sigma  ^{2} (u(t_{j}(\delta),x))\right]\\	&&\mathbf{E}\left[\left( N ^{-\frac{1}{\alpha}}(\widetilde{\Delta} _{1} u_{0} (t_{i}, \delta))^{2}- C_{0, \alpha}^{2} \frac{1}{N}\right)\left( N ^{-\frac{1}{\alpha}}(\widetilde{\Delta} _{1} u_{0} (t_{j}, \delta))^{2}- C_{0, \alpha}^{2} \frac{1}{N}\right)\right] \\
		&\leq & C \sum_{i,j=0; i-j\geq [N ^{-\beta +1}]+2}^{N-1} \left| \mathbf{E}\left[\left( N ^{-\frac{1}{\alpha}}(\widetilde{\Delta} _{1} u_{0} (t_{i}, \delta))^{2}- C_{0, \alpha}^{2} \frac{1}{N}\right)\right. \right.\\
		&&\times \left. \left. \left( N ^{-\frac{1}{\alpha}}(\widetilde{\Delta} _{1} u_{0} (t_{j}, \delta))^{2}- C_{0, \alpha}^{2} \frac{1}{N}\right)\right] \right|.
	\end{eqnarray*}
	Now, we use the estimate (\ref{14s-1})  in Lemma \ref{ll2} to get, with $\rho_{H}$ given by (\ref{ro})
	\begin{eqnarray}
		l_{2,N} ^{(2,2)}&\leq &C\frac{1}{N^{2}} \left( \sum_{i,j=1; i-j\geq [N ^{-\beta +1}]+2}^{N-1} \left( \rho ^{2}_{ \frac{\alpha -1}{2\alpha}}(i-j)+ i ^{H-1} \right)+1\right)\nonumber\\
		&\leq & C \frac{1}{N^{2}}\left( \sum_{i,j=1; i-j\geq1}^{N-1} \left( \rho ^{2}_{ \frac{\alpha -1}{2\alpha}}(i-j)+ i ^{H-1} \right)+1\right)\nonumber \\
		&\leq &  C \left( \frac{1}{N ^{2}} \sum _{j=1} ^{N} \sum _{i=j+1} ^{N} \rho ^{2}_{ \frac{\alpha -1}{2\alpha}}(i-j) + \frac{1}{N} \sum_{i=1} ^{N-1} i ^{H-1} +\frac{1}{ N ^{2}}\right)\nonumber \\
		&=& C\left( \frac{1}{N ^{2}} \sum _{j=1} ^{N} \sum_{k=1} ^{N-1-j} \rho ^{2}_{ \frac{\alpha -1}{2\alpha}}(k) +\frac{1}{N} \sum_{i=1} ^{N-1} i ^{H-1} +\frac{1}{ N ^{2}}\right)\nonumber\\
		&\leq & C \left( \frac{1}{N} \sum_{i=1}^{N-1}\rho^{2} _{\frac{\alpha-1}{2\alpha}}(i) + \frac{1}{N} \sum_{i=1} ^{N-1} i ^{H-1} +\frac{1}{ N ^{2}}\right)\nonumber\\
		&\leq & C\left( N ^{-1} + N ^{-\frac{\alpha +1}{2\alpha}}\right) \leq C N ^{-\frac{\alpha +1}{2\alpha}}.	\label{27i-1}	 
	\end{eqnarray}
	We have obtained the bound
	\begin{equation}
		\label{3f-2}	\mathbf{E} \vert L_{2,N}\vert ^{2} \leq C N ^{\frac{-\alpha-1}{{2\alpha}}}, \mbox{ so } \mathbf{E} \vert L_{2,N}\vert \leq C N ^{\frac{-\alpha-1}{{4\alpha}}},
	\end{equation}
	
	\noindent	Finally, the term $ L_{3,N}$ can be estimated in the following way. 
	\begin{eqnarray*}
		L_{3,N} &=& C_{0, \alpha} ^{2} \sum_{i=0} ^{N-1} \int_{t_{i}}^{t_{i+1}} \left( \sigma ^{2} (u(t_{i}(\delta), x))-  \sigma ^{2} (u(s, x))\right) ds 
	\end{eqnarray*}
	and then, by (\ref{2i-2}) and (\ref{td}),
	\begin{eqnarray}
		\mathbf{E} \vert L_{3, N}\vert &\leq & C \sum_{i=0} ^{N-1} \int_{t_{i}}^{t_{i+1}} \left( 	\mathbf{E} \left( \sigma (u(t_{i}(\delta), x))-\sigma (u(s,x))\right) ^{2} \right) ^{\frac{1}{2}}ds\nonumber\\
		&\leq & C  \sum_{i=0} ^{N-1} \int_{t_{i}}^{t_{i+1}} \left( 	\mathbf{E} \vert u(t_{i}( \delta), x)-u(s,x) \vert ^{2} \right) ^{\frac{1}{2}}ds\nonumber\\
		&\leq & C   \sum_{i=0} ^{N-1} \int_{t_{i}}^{t_{i+1}}\vert t_{i}(\delta)-s\vert ^{\frac{\alpha -1}{2\alpha}}ds\\
		&\leq & C \left( \frac{1}{N^{\beta}} \right) ^{\frac{\alpha-1}{2\alpha}}=C \left( \frac{1}{N}\right) ^{\frac{ \alpha-1}{2\alpha+1}}\to _{N \to \infty} 0. \label{3f-12}
	\end{eqnarray}
	By combining the estimates (\ref{3f-11}), (\ref{3f-2}) and (\ref{3f-12}) we notice that the dominant term is  $ L_{1,N} $ and we  obtain the inequality (\ref{3f-5}). \qed
	
	\subsection{Higher order variations}
	Let $u_{0}$  be given by (\ref{u0}), i.e. $u_{0}$ is the solution to the fractional stochastic heat equation with additive spce-time white noise (that is,  $\sigma (x)=1$ for every $x \in \mathbb{R}$).
	Recall that $u_{0}$ is a perturbed fBm with Hurst parameter $ H= \frac{\alpha-1}{2\alpha}$ (see the decomposition  (\ref{deco})). The purpose is to analyze the $\frac{1}{H}$-temporal variation of the solution to the heat equation defined by (\ref{mild1}). 
	
	We recall the following result (see \cite{Mahtu2} or \cite{T2}):
	\begin{equation}\label{boa}
		U_{N,x} (u_{0}) =\sum_{i=0} ^{N-1} \left| u_{0} (t_{i+1}, x)-u_{0} (t_{i}, x) \right| ^{\frac{2\alpha}{\alpha-1}}\to _{N \to \infty} B_{0, \alpha}:= C_{0, \alpha} ^{\frac{2\alpha}{\alpha-1}}	\mathbf{E} \vert Z \vert ^{\frac{2\alpha}{\alpha-1}} \mbox{ in } L ^{1} (\Omega),
	\end{equation}
	where $Z \sim N(0,1)$ and $ C_{0, \alpha}$ is the strictly positive constant appearing in the decomposition (\ref{deco}). 
	
	The goal  is to prove a similar result for the solution $u$ to the nonlinear fractional stochastic heat equation (\ref{2}). We set, for every $N\geq 1$, 
	
	\begin{equation}
		\label{unu}
		U_{N,x}(u)= \sum_{i=0} ^{N-1} \left| u (t_{i+1}, x)-u (t_{i}, x) \right| ^{\frac{2\alpha}{\alpha-1}}.
	\end{equation}
	Below we analyze the limit behavior the sequence (\ref{unu}). We assume in the sequel that
	\begin{equation}\label{28i-1}
		\frac{2\alpha}{\alpha-1} \mbox{ is an  even integer. }
	\end{equation}

	We need need another technical lemma concerning the increment of the perturbed fBm.

	\begin{lemma}\label{ll3}
		Let $(X_{t}, t\geq 0)$ be a perturbed fBm defined by (\ref{pfBm}) and let $t_{i}, i-0,1,..,N$ be given by (\ref{ti}). Fix $N\geq 1$ and $i=1,...,N$.  Then
		\begin{enumerate}

			\item If  $p=\frac{1}{H}$ is an integer,
			\begin{equation}\label{17s-6}
				\mathbf{E} \left( 	  \left( X_{t_{i+1}}-X_{t_{i}}\right) ^{p} -\frac{A_{0,p} }{N} \right) ^{2} \leq C \frac{1}{N ^{2}},
			\end{equation}
			where 
			\begin{equation}
				\label{bop}
				A_{0,p}= C_{0} ^{p} 	\mathbf{E} (Z^{p}), \mbox{ where }Z \sim N(0,1).
			\end{equation} 
			\item If $p=\frac{1}{H}$ is an  integer, then for every $j=1,..., N$ with $j\not=i$, 
			\begin{eqnarray}
				&&	\left| 	\mathbf{E}\left(  \left( X_{t_{i+1}}-X_{t_{i}}\right) ^{p} -\frac{A_{0,p} }{N} \right) \left(  \left( X_{t_{j+1}}-X_{t_{j}}\right) ^{p} -\frac{A_{0,p} }{N} \right) \right| \nonumber\\
				&&\leq C\frac{1}{N ^{2}}(\rho_{H}(i-j) + i ^{H-1}+ j^{H-1}),\label{17s-2}
			\end{eqnarray}
			with $\rho_{H}$ defined by (\ref{ro}).
		\end{enumerate}
	\end{lemma}
	{\bf Proof: } We use the decomposition 
	\begin{eqnarray}\label{8f-4}
		\left( X_{t_{i+1}}-X_{t_{i}}\right) ^{p} -\frac{ A_{0, p}}{N}= C_{0}^{p}\left( B ^{H}_{t_{i+1}}-B ^{H} _{t_{i}}\right) ^{p} -\frac{ B_{0, p}}{N} + R_{i, N},
	\end{eqnarray}
	where 
	\begin{eqnarray*}
		R_{i,N} =\sum_{k=1} ^{p-1}C_{p} ^{k}C_{0}^{p-k} \left( B^{H}_{t_{i+1}}-B ^{H}_{t_{i}}\right) ^{p-k} \left( Y_{t_{i+1}}-Y_{t_{i}}\right) ^{k} + \left( Y_{t_{i+1}}-Y_{t_{i}}\right) ^{p} 
	\end{eqnarray*}
	We have  the bound
	\begin{eqnarray}
		&&	\mathbf{E}R_{i, N} ^{2}\nonumber\\
		&&\leq \sum_{k=1}^{p-1} c_{p} \left( \mathbf{E}\left( B^{H}_{t_{i+1}}-B ^{H}_{t_{i}}\right) ^{4p-4k}\right) ^{\frac{1}{2}}\left( \mathbf{E}\left( Y_{t_{i+1}}- Y_{t_{i}}\right) ^{4k}\right) ^{\frac{1}{2}} + \mathbf{E} \left( Y_{t_{i+1}}-Y_{t_{i}}\right) ^{2p}\nonumber  \\
		&&
		\leq C\frac{ i^{2H-2}}{N ^{2pH}}=C\frac{ i^{2H-2}}{N ^{2}}\label{17s-4}
	\end{eqnarray}
	and  we also notice that,  if $Z$ stands for a $N(0,1)$ random variable, 
	\begin{eqnarray}
		\mathbf{E} \left( C_{0} ^{p}	\left( B ^{H}_{t_{i+1}}-B ^{H} _{t_{i}}\right) ^{p} -\frac{ A_{0, p}}{N} \right) ^{2} =C_{0} ^{2p} \frac{1}{ N ^{2}}	\mathbf{E} \vert Z ^{p}-	\mathbf{E}Z^{p} \vert ^{2} = C \frac{1}{N ^{2}}.\label{17s-5}
	\end{eqnarray}
	By (\ref{17s-4}) and (\ref{17s-5}), we obtain the inequality (\ref{17s-6}).   For the bound (\ref{17s-2}), we first prove that 
	\begin{eqnarray}
		&&	\left| 	\mathbf{E}\left(  \left( B ^{H}_{t_{i+1}}-B ^{H}_{t_{i}}\right) ^{p} -\frac{\mathbf{E} (Z ^{p}) }{N} \right) \left(  \left( B ^{H}_{t_{j+1}}-B ^{H} _{t_{j}}\right) ^{p} -\frac{\mathbf{E}(Z ^{p}) }{N} \right) \right| \nonumber\\
		&&\leq C\frac{1}{N ^{2}}\rho_{H}(i-j),	\label{17s-3}
	\end{eqnarray}
	where $Z$ is a standard Gaussian random variable. To do it, we  observe that, by the scaling property of the fBm,  
	\begin{eqnarray*}
		&&	\mathbf{E}\left(  \left( B ^{H}_{t_{i+1}}-B ^{H}_{t_{i}}\right) ^{p} -\frac{\mathbf{E} (Z ^{p}) }{N} \right) \left(  \left( B ^{H}_{t_{j+1}}-B ^{H} _{t_{j}}\right) ^{p} -\frac{\mathbf{E}(Z ^{p}) }{N} \right)\\
		&&= \frac{1}{N ^{2}}\mathbf{E}\left[  ( (\Delta B^{H}_{i}) ^{p}-\mathbf{E} (Z^{p}) )((\Delta B^{H}_{j}) ^{p} -\mathbf{E}(Z ^{p}))\right],
	\end{eqnarray*}
	with $\Delta B^{H}_{i}= B^{H}_{i+1}- B ^{H}_{i}$.  By Gaussian regression, we can write
	\begin{equation*}
		\Delta B^{H}_{i}=\rho_{H}(i-j) \Delta B^{H}_{j}+ \sqrt{1-\rho_{H} ^{2}(i-j)}G,
	\end{equation*}
	with $ G\sim N(0,1)$ independent of $\Delta B^{H}_{j}$. This leads to 
	\begin{eqnarray*}
		&&\mathbf{E}\left[  ( (\Delta B^{H}_{i}) ^{p}-\mathbf{E} (Z^{p}) )((\Delta B^{H}_{j}) ^{p} -\mathbf{E}(Z ^{p}))\right]\\
		&=&\mathbf{E}\left[ \left(  \left( \rho_{H}(i-j) \Delta B^{H}_{j}+ \sqrt{1-\rho_{H} ^{2}(i-j)}G\right) ^{p} -\mathbf{E}(Z ^{p})\right)((\Delta B^{H}_{j}) ^{p} -\mathbf{E}(Z ^{p}))\right] \\
		&=& \mathbf{E}\left[  \left( \sum_{k=1} ^{p} C_{p}^{k} \rho_{H}(i-j) ^{k} (\Delta B^{H}_{j}) ^{k} \sqrt{1-\rho_{H} ^{2}(i-j)}^{p-k} G ^{p-k} -\mathbf{E}(Z ^{p})\right) ((\Delta B^{H}_{j}) ^{p} -\mathbf{E}(Z ^{p}))\right]\\
		&=& \sum_{k=0}^{p}C_{p} ^{k}  \rho_{H}(i-j) ^{k} (1-\rho_{H} ^{2}(i-j))^{\frac{p-k}{2}}\left( \mathbf{E} (\Delta B ^{H}_{j} ) ^{p+k} G ^{p-k} -\mathbf{E} \left( (\Delta B^{H}_{j})^{k}G ^{p-k}\right) \mathbf{E} (Z ^{p}) \right), 
	\end{eqnarray*}
	and, since $ \Delta B ^{H}_{j} $ and $ G$ are independent standard normal random variables,

	\begin{eqnarray*}
		&&\mathbf{E}\left[  ( (\Delta B^{H}_{i}) ^{p}-\mathbf{E} (Z^{p}) )((\Delta B^{H}_{j}) ^{p} -\mathbf{E}(Z ^{p}))\right]\\
		&=&\sum_{k=1}^{p}C_{p} ^{k}  \rho_{H}(i-j) ^{k} (1-\rho_{H} ^{2}(i-j))^{\frac{p-k}{2}} \mathbf{E} ( Z ^{p-k}) \left( \mathbf{E} (Z ^{p+k}) - \mathbf{E} ( Z ^{p}) \mathbf{E} (Z ^{k} )\right)\\
		&=&\sum_{k=1}^{p}C_{p} ^{k}  \rho_{H}(i-j) ^{k} (1-\rho_{H} ^{2}(i-j))^{\frac{p-k}{2}} \mathbf{E} ( Z ^{p-k}) \left( \mathbf{E} (Z ^{p+k} )- \mathbf{E} ( Z ^{p}) \mathbf{E} (Z ^{k})\right)\\
		&\leq & C \rho_{H} (i-j),
	\end{eqnarray*}
	where we used   $\vert \rho_{H}(i-j)\vert \leq 1$ for all $i, j$. So (\ref{17s-3}) follows. The proof of (\ref{17s-2}) can be then deduced from the decomposition (\ref{8f-4}) and the relations  (\ref{23i-1}), (\ref{17s-3}).
	\qed
	
	We state and prove our result concerning the limit behavior of the sequence (\ref{unu}). 
	
	\begin{theorem} \label{tt2} Consider the sequence $ (U_{N,x}(u), N\geq 1)$ given by (\ref{unu}).  Then
		\begin{equation}\label{13f-1}
			\mathbf{E} \left| U_{N,x}(u)- B_{0, \alpha}\int_{0} ^{1} \sigma ^{\frac{ 2\alpha}{\alpha-1}} (u(s,x))ds\right| \leq C\left(  \frac{1}{N}\right) ^{ \frac{ (\alpha-1) (2\alpha-1)}{2\alpha (2\alpha+1)}}.
		\end{equation}
		In particular,  $ (U_{N,x}(u), N\geq 1)$ converges to $ B_{0, \alpha}\int_{0} ^{1} \sigma ^{\frac{ 2\alpha}{\alpha-1}} (u(s,x))ds$ in $ L ^{1} (\Omega)$ as $N \to \infty$. The constant $ B_{0, \alpha}$ is defined in (\ref{boa}). 
		
	\end{theorem}
	{\bf Proof: } We use a decomposition similar to that used in the proof of Theorem \ref{tt1}. We write 
	\begin{eqnarray}
		&&	U_{N,x}(u)- B_{0, \alpha}\int_{0} ^{1} \sigma ^{\frac{ 2\alpha}{\alpha-1}} (u(s,x))ds\nonumber\\
		&=& \sum_{i=0} ^{N-1} \left( \left| u (t_{i+1}, x)-u (t_{i}, x) \right| ^{\frac{2\alpha}{\alpha-1}}-\sigma ^{\frac{2\alpha}{\alpha -1}}(u(t_{i}(\delta),x))\left( \widetilde{\Delta}_{1} u_{0} (t_{i}, \delta) \right) ^{\frac{2\alpha}{\alpha -1}}\right)\nonumber\\
		&&+\sum_{i=0}^{N-1} \sigma ^{\frac{2\alpha}{\alpha -1}}(u(t_{i}(\delta),x)) \left[ \left( \widetilde{\Delta}_{1} u_{0} (t_{i}, \delta) \right) ^{\frac{2\alpha}{\alpha -1}}-\frac{ B_{0, \alpha}}{N} \right]\nonumber\\
		&&+ B_{0, \alpha} \left[ \frac{1}{N}\sum _{i=0} ^{N-1}  \sigma ^{\frac{2\alpha}{\alpha -1}}(u(t_{i}(\delta),x))-\int_{0} ^{1} \sigma ^{\frac{ 2\alpha}{\alpha-1}} (u(s,x))ds\right]\nonumber\\
		&=:& A_{1, N}+ A_{2, N}+ A_{3, N}. \label{deco1}
	\end{eqnarray}
	We analize separately the  behavior of the three summands from above as $N \to \infty$.  Let $p=\frac{2\alpha}{\alpha-1}$ and recall that $p$ is an  even integer ($p\geq 2$). We have
	\begin{equation*}
		A_{1,N}= \sum_{i=0} ^{N-1} \left[  \left( u (t_{i+1}, x)-u (t_{i}, x) \right)-\sigma (u(t_{i}(\delta))) \widetilde{\Delta}_{1} u_{0} (t_{i}, \delta) \right] R_{i, N}
	\end{equation*}
	with
	\begin{equation*}
		R_{i, N}=\sum_{k=1} ^{p-1} \left( u (t_{i+1}, x)-u (t_{i}, x) \right)^{p-1-k}\sigma^{k} (u(t_{i}(\delta)))\left( \widetilde{\Delta}_{1} u_{0} (t_{i}, \delta)\right) ^{k}.
	\end{equation*}
	From (\ref{2i-2}) and (\ref{8f-1}), it  holds that, for $N \geq 1$ and for $i=1,...,N$, 
	\begin{equation*}
		\mathbf{E} R_{i, N}^{2} \leq C \frac{1}{N ^{2(p-1)H}}=N ^{2H-2},
	\end{equation*}
	where $ H=\frac{\alpha -1}{2\alpha}$.  Thus, by Proposition \ref{pp4},
	\begin{eqnarray*}
		\mathbf{E} \vert A_{1, N}\vert &\leq & \sum_{i=0} ^{N-1} \left( 	\mathbf{E} \left[  \left( u (t_{i+1}, x)-u (t_{i}, x) \right)-\sigma (u(t_{i}(\delta))) \widetilde{\Delta}_{1} u_{0} (t_{i}, \delta) \right] ^{2} \right) ^{\frac{1}{2}} (	\mathbf{E}R_{i, N}^{2}) ^{\frac{1}{2}}\\
		&\leq & C\sum_{i=0} ^{N-1}\left( \frac{1}{N} \right)^{\frac{ 2(\alpha-1)}{2\alpha+1}} N ^{H-1}\leq C\left(  \frac{1}{N}\right) ^{ \frac{ (\alpha-1) (2\alpha-1)}{2\alpha (2\alpha+1)}}.
	\end{eqnarray*}
	Concerning $A_{3, N}$,
	\begin{eqnarray*}
		A_{3,N} &=& B_{0, \alpha}\left( \sum_{i=0} ^{N-1} \int_{t_{i}}^{t_{i+1}} \left(  \sigma ^{\frac{2\alpha}{\alpha -1}}(u(t_{i}(\delta),x))- \sigma ^{\frac{2\alpha}{\alpha -1}}(u(s,x))\right) ds \right) \\&=& B_{0, \alpha} \left( \sum_{i=0}^{N-1}  \int_{t_{i}}^{t_{i+1}} ds \left(  \sigma (u(t_{i}(\delta),x))- \sigma (u(s,x))\right)  R_{i, N}(s) \right),
	\end{eqnarray*}
	with 
	\begin{equation*}
		R_{i, N}(s)=\sum_{k=1} ^{p-1}  \sigma ^{p-1-k}(u(t_{i}(\delta),x)) \sigma ^{k}(u(s,x)),
	\end{equation*}
	and clearly, by (\ref{2i-3}) and (\ref{3i-5}),
	\begin{equation*}
		\mathbf{E} R_{i, N}(s) ^{2} \leq C,
	\end{equation*}
	for every $N\geq 1, i=1,...,N$ and $s\in (0, 1)$.  Therefore,
	\begin{eqnarray*}
		\mathbf{E} \vert A_{3, N}\vert &\leq & C  \sum_{i=0} ^{N-1} \int_{t_{i}}^{t_{i+1}}ds \left( 	\mathbf{E} \left(  \sigma (u(t_{i}(\delta),x))- \sigma (u(s,x))\right) ^{2} \right) ^{\frac{1}{2}} (	\mathbf{E} R_{i,N}(s))^{\frac{1}{2}}\\
		&\leq & C  \sum_{i=0} ^{N-1} \int_{t_{i}}^{t_{i+1}}ds \left( 	\mathbf{E} \left(  \sigma (u(t_{i}(\delta),x))- \sigma (u(s,x))\right) ^{2} \right) ^{\frac{1}{2}}\\
		&\leq &  C  \sum_{i=0} ^{N-1} \int_{t_{i}}^{t_{i+1}}ds\left( 	\mathbf{E} \left(  u(t_{i}(\delta), x)-  u(s,x)\right) ^{2} \right) ^{\frac{1}{2}}\\
		&\leq &  C  \sum_{i=0} ^{N-1} \int_{t_{i}}^{t_{i+1}}ds \vert t_{i}(\delta) -s\vert ^{\frac{ \alpha -1}{2\alpha}}\leq C \left( \frac{1}{N^{\beta}} \right) ^{\frac{ \alpha -1}{2\alpha}}\leq C \left( \frac{1}{N}\right) ^{\frac{\alpha-1}{2\alpha+1}}.
	\end{eqnarray*}
	Let us finally deal with the summand $ A_{2, N}$. We estimate its $L^{2}(\Omega)$-norm. We have 
	\begin{eqnarray*}
		\mathbf{E} A_{2, N} ^{2} &=&	\mathbf{E} \sum_{i,j=0}^{N-1} \sigma ^{\frac{2\alpha}{\alpha -1}}(u(t_{i}(\delta),x)) \sigma ^{\frac{2\alpha}{\alpha -1}}(u(t_{j}(\delta),x)) \\
		&&\left[ \left( \widetilde{\Delta}_{1} u_{0} (t_{i}, \delta) \right) ^{\frac{2\alpha}{\alpha -1}}-\frac{ B_{0, \alpha}}{N} \right] \left[ \left( \widetilde{\Delta}_{1} u_{0} (t_{j}, \delta) \right) ^{\frac{2\alpha}{\alpha -1}}-\frac{ B_{0, \alpha}}{N} \right]\\
		&=& 	\mathbf{E} \sum_{i=0}^{N-1} \sigma ^{\frac{4\alpha}{\alpha -1}}(u(t_{i}(\delta),x))\left[ \left( \widetilde{\Delta}_{1} u_{0} (t_{i}, \delta) \right) ^{\frac{2\alpha}{\alpha -1}}-\frac{ B_{0, \alpha}}{N} \right] ^{2}\\
		&&+ 2	\mathbf{E} \sum_{i,j=0; i>j}^{N-1} \sigma ^{\frac{2\alpha}{\alpha -1}}(u(t_{i}(\delta),x)) \sigma ^{\frac{2\alpha}{\alpha -1}}(u(t_{j}(\delta),x)) \\
		&&\left[ \left( \widetilde{\Delta}_{1} u_{0} (t_{i}, \delta) \right) ^{\frac{2\alpha}{\alpha -1}}-\frac{ B_{0, \alpha}}{N} \right] \left[ \left( \widetilde{\Delta}_{1} u_{0} (t_{j}, \delta) \right) ^{\frac{2\alpha}{\alpha -1}}-\frac{ B_{0, \alpha}}{N} \right]\\
		&=:& a_{2, N} ^{(1)} + a_{2, N}^{(2)}.
	\end{eqnarray*}
	We use the independence of $\widetilde{\Delta}_{1} u_{0} (t_{j}, \delta)$ and  $\sigma ^{\frac{2\alpha}{\alpha -1}}(u(t_{j}(\delta),x)),$ the inequality  (\ref{3i-5})  and Lemma \ref{ll3} to get 
	\begin{eqnarray*}
		a_{2, N} ^{(1)} &= & \sum_{i=0}^{N-1} 	\mathbf{E} \sigma ^{\frac{4\alpha}{\alpha -1}}(u(t_{i}(\delta),x))	\mathbf{E}\left[ \left( \widetilde{\Delta}_{1} u_{0} (t_{i}, \delta) \right) ^{\frac{2\alpha}{\alpha -1}}-\frac{ B_{0, \alpha}}{N} \right] ^{2}\\
		&\leq & C 	\sum_{i=0}^{N-1} 	\mathbf{E}\left[ \left( \widetilde{\Delta}_{1} u_{0} (t_{i}, \delta) \right) ^{\frac{2\alpha}{\alpha -1}}-\frac{ B_{0, \alpha}}{N} \right] ^{2}\\
		&\leq & C \sum_{i=0}^{N-1}\frac{1}{N ^{2}}=C \frac{1}{N}.
	\end{eqnarray*}
	For the summand $a_{2, N} ^{(2)}$, we divide it into two sums, as in the proof of Theorem 1. A first sum will contain all the terms such that $i-j\leq [ N ^{-\beta +1}]+1$. This first sum will be denoted by $a_{2, N} ^{(2,1)}$. In the second sum we include all the terms with $i-j\geq [ N ^{-\beta +1}]+2$. This second sum will be denoted by $a_{2, N} ^{(2,2)}$. 
	
	We treat $ a_{2, N} ^{(2,1)} $ in the following way.
	\begin{eqnarray*}
		a_{2, N} ^{(2,1)} &\leq & \sum_{i,j=0; i-j\leq [N ^{-\beta+1}]+1}^{N-1} \left( 	\mathbf{E} \sigma ^{\frac{4\alpha}{\alpha -1}}(u(t_{i}(\delta),x)) \left[ \left( \widetilde{\Delta}_{1} u_{0} (t_{i}, \delta) \right) ^{\frac{2\alpha}{\alpha -1}}-\frac{ B_{0, \alpha}}{N} \right] ^{2}\right)^{\frac{1}{2}}\\
		&&\times   \left( 	\mathbf{E} \sigma ^{\frac{4\alpha}{\alpha -1}}(u(t_{j}(\delta),x)) \left[ \left( \widetilde{\Delta}_{1} u_{0} (t_{j}, \delta) \right) ^{\frac{2\alpha}{\alpha -1}}-\frac{ B_{0, \alpha}}{N} \right] ^{2}\right)^{\frac{1}{2}}\\
		&=& \sum_{i,j=0; i-j\leq [N ^{-\beta+1}]+1}^{N-1} \left( 	\mathbf{E} \sigma ^{\frac{4\alpha}{\alpha -1}}(u(t_{i}(\delta),x)) 	\mathbf{E}\left[ \left( \widetilde{\Delta}_{1} u_{0} (t_{i}, \delta) \right) ^{\frac{2\alpha}{\alpha -1}}-\frac{ B_{0, \alpha}}{N} \right] ^{2}\right)^{\frac{1}{2}}\\
		&&\times   \left( 	\mathbf{E} \sigma ^{\frac{4\alpha}{\alpha -1}}(u(t_{j}(\delta),x)) 	\mathbf{E}\left[ \left( \widetilde{\Delta}_{1} u_{0} (t_{j}, \delta) \right) ^{\frac{2\alpha}{\alpha -1}}-\frac{ B_{0, \alpha}}{N} \right] ^{2}\right)^{\frac{1}{2}}.
	\end{eqnarray*}
	The bound (\ref{3i-5}) and  Lemma \ref{ll3}  point 1. lead to the estimate
	\begin{eqnarray*}
		a_{2, N} ^{(2,1)} &\leq &C \sum_{i,j=0; i-j\leq [N ^{-\beta+1}]+1}^{N-1} \left(  	\mathbf{E}\left[ \left( \widetilde{\Delta}_{1} u_{0} (t_{i}, \delta) \right) ^{\frac{2\alpha}{\alpha -1}}-\frac{ B_{0, \alpha}}{N} \right] ^{2}\right)^{\frac{1}{2}}\\
		&&\times   \left(  	\mathbf{E}\left[ \left( \widetilde{\Delta}_{1} u_{0} (t_{j}, \delta) \right) ^{\frac{2\alpha}{\alpha -1}}-\frac{ B_{0, \alpha}}{N} \right] ^{2}\right)^{\frac{1}{2}}\\
		&\leq & C \sum_{i,j=0; i-j\leq [N ^{-\beta+1}]+1}^{N-1} \frac{1}{N ^{2}} \leq C N ^{-\beta},
	\end{eqnarray*}
	with $\beta$ from (\ref{td}).  Concerning the summand $ a_{2, N} ^{(2,2)}$, we have since $\widetilde{\Delta}_{1} u_{0} (t_{i}, x), \widetilde{\Delta}_{1} u_{0} (t_{j}, \delta)$ are independent of $ (\sigma(u(t_{i}(\delta),x)), \sigma ^{\frac{2\alpha}{\alpha -1}}(u(t_{j}(\delta),x)))$,
	\begin{eqnarray*}
		&&	a_{2, N} ^{(2,2)}= \sum_{i,j=0; i-j\geq [N ^{-\beta+1}]+2}^{N-1}\\
		&& \left|  	\mathbf{E} \sigma ^{\frac{2\alpha}{\alpha -1}}(u(t_{i}(\delta),x)) \sigma ^{\frac{2\alpha}{\alpha -1}}(u(t_{j}(\delta),x)) 
		\right| \\
		&&	\left| 	\mathbf{E}  \left[ \left( \widetilde{\Delta}_{1} u_{0} (t_{j}, \delta) \right) ^{\frac{2\alpha}{\alpha -1}}-\frac{ B_{0, \alpha}}{N} \right]	\left[ \left( \widetilde{\Delta}_{1} u_{0} (t_{j}, \delta) \right) ^{\frac{2\alpha}{\alpha -1}}-\frac{ B_{0, \alpha}}{N} \right]\right| \\
		&\leq &C \sum_{i,j=0; i-j\geq [N ^{-\beta+1}]+2}^{N-1} 	\left| 	\mathbf{E}  \left[ \left( \widetilde{\Delta}_{1} u_{0} (t_{j}, \delta) \right) ^{\frac{2\alpha}{\alpha -1}}-\frac{ B_{0, \alpha}}{N} \right]	\left[ \left( \widetilde{\Delta}_{1} u_{0} (t_{j}, \delta) \right) ^{\frac{2\alpha}{\alpha -1}}-\frac{ B_{0, \alpha}}{N} \right]\right| 
	\end{eqnarray*}
	and via (\ref{17s-2}) in Lemma \ref{ll3}, as we proceeded in (\ref{27i-1})
	\begin{eqnarray*}
		a_{2, N} ^{(2,2)}&\leq &  C \frac{1}{N^{2}}\left(1+ \sum_{i,j=1; i-j\geq1}^{N-1} \left( \rho _{ \frac{\alpha -1}{2\alpha}}(i-j)+ i ^{H-1} \right)\right)\\
		&\leq & C \left( \frac{1}{ N ^{2}}+\frac{1}{N} \sum_{i=1}^{N-1}\rho _{\frac{\alpha-1}{2\alpha}}(i) + \frac{1}{N} \sum_{i=1} ^{N-1} i ^{H-1} \right)\\
		&\leq & C\left( N ^{-1} + N ^{-\frac{\alpha +1}{2\alpha}}\right) \leq C N ^{-\frac{\alpha +1}{2\alpha}}.		 
	\end{eqnarray*}
	By comparing the estimates for $ 	\mathbf{E}\vert A_{1,N}\vert , 	\mathbf{E}\vert A_{3, N}\vert $ and $ 	\mathbf{E} \vert A_{2, N}\vert$ (which is less than $ (	\mathbf{E} \vert A_{2, N}\vert ^{2} )^{\frac{1}{2}}$), we obtain (\ref{13f-1}). In particular, we notice that this biggest term in the decomposition (\ref{deco1}) (in the sense of the $ L ^{1}(\Omega)$-norm) is $ A_{1, N}$. 
	\qed

	\begin{remark}\label{rem1}
		\begin{enumerate}
			
			\item We impose the above condition (\ref{28i-1})  because our proofs are based on the behaviour of the power variations of the fBm where the study has been done for integer powers  (see e.g. the survey of results  in the introduction of \cite{NNT}).  The condition that $\frac{2\alpha}{\alpha-1}$ is even is needed for technical reasons (for instance, this allows to treat the summand $A_{1,N}$ in the proof of Theorem \ref{tt2} via the estimate obtained in Proposition \ref{pp4}).

			\item 
			The result in Theorem \ref{tt2}, as well as the one in Theorem \ref{tt1}, can be easily extended to any time interval $ [A_{1}, A_{2}]$ (instead of $[0, 1]$). Indeed, if one considers the equidistant partition of the interval $[A_{1}, A_{2}]$,
			\begin{equation*}
				t_{i}= A_{1}+ \frac{i}{N} (A_{2}-A_{1})
			\end{equation*}
			one has, in $ L ^{1} (\Omega)$, 
			\begin{equation}
				\label{13f-2}
				N ^{-\frac{1}{\alpha}} \sum_{i=0} ^{N-1} \left( u(t_{i+1}, x)-u(t_{i}, x)\right) ^{2}\to_{N \to \infty} C_{0, \alpha } ^{2}(A_{2}- A_{1}) ^{2H-1}\int_{A_{1}}^{A_{2}} \sigma ^{2} (u(s, x))ds
			\end{equation}
			and, if $\frac{2\alpha}{\alpha- 1}$ is an even integer,
			\begin{equation}\label{13f-3}
				\sum_{i=0} ^{N-1} \left| u (t_{i+1}, x)-u (t_{i}, x) \right| ^{\frac{2\alpha}{\alpha-1}}\to_{N \to \infty}B_{0, \alpha}\int_{A_{1}} ^{A_{2}} \sigma ^{\frac{ 2\alpha}{\alpha-1}} (u(s,x))ds,
			\end{equation}
			with $C_{0, \alpha}, B_{0, \alpha}$ given by (\ref{deco}) and (\ref{boa}), respectively. 
			
			\item If $\alpha =2$ (i.e. $-(-\Delta) ^{\frac{\alpha}{2}}$ is the standard Laplacian operator), then it follows from Theorem \ref{tt2} that the solution to the (\ref{2}) has a nontrivial quartic variation. This corresponds to some known facts from the literature, see \cite{PoTr} or \cite{Sw}.   In particular, the right-hand side of the bound (\ref{13f-1}) reads $C N ^{-\frac{3}{20}}$ (this corresponds with the bound obtained in the proof of Proposition 5.3 in \cite{PoTr}). 
			
			\item By comparing the bounds (\ref{3f-5}) and (\ref{13f-1}), we notice  that the quadratic variation (\ref{vnxu}) and the power variation (\ref{unu}) have the same rate of convergence in $ L ^{1} (\Omega)$. 
		\end{enumerate} 
	\end{remark}

	\section{Application to parameter estimation}
	We now apply our main results stated in Theorem \ref{tt1} and Theorem \ref{tt2} to the estimation of the  parameters that appear in the model (\ref{intro-1}). In the first part, we illustrate how the anomality parameter $\alpha$ (which appear in (\ref{intro-1}) with the fractional Laplacian operator) may be identified from the observation of the solution to (\ref{intro-1}) at discrete times and at a fixed spatial point.  In a second part, we derive a consistent estimator for the drift parameter $\theta$.

	\subsection{The variations of the parametrized heat equation}
	In this part, we consider the parametrized stochastic heat equation 
	\begin{equation}
		\label{1}
		\frac{\partial u_{\theta }}{\partial t}(t,x)= -\theta (-\Delta) ^{ \frac{\alpha}{2}}u_{\theta}(t,x)+  \sigma (u_{\theta }(t,x)) \dot{W} (t,x), \hskip0.4cm t\geq 0, x\in \mathbb{R},
	\end{equation}
	with vanishing initial value $u_{\theta }(0, x)=0$ for every $x\in \mathbb{R}$. As above, $W$ is a space-time white noise  (a centered Gaussian random field characterized by its covariance (\ref{covW})) and we assume that $\theta >0$ and the parameter $\alpha $ belongs to the interval $(1,2]$. 
	
	We start by deducing  from Theorem \ref{tt1} and  Theorem \ref{tt2} the behavior of the quadratic and higher-order variations for solution to the parametrized heat equation (\ref{1}).  To this end, we use a standard procedure. We define the transform
	\begin{equation*}
		v_{\theta }(t,x)= u_{\theta }\left( \frac{t}{\theta}, x\right) \mbox{ for every } t\geq 0, x\in \mathbb{R},
	\end{equation*}
	and  it follows that $v_{\theta} $ satisfies the equation (see \cite{PoTr}, \cite{GT} or \cite{T2})
	\begin{equation*}
		v_{\theta } (t,x)= \int_{0} ^{t} \int_{\mathbb{R}} G_{\alpha} (t-s, x-y)\widehat{\sigma} (v_{\theta}(s,y))\widetilde{W} (ds, dy),
	\end{equation*}
	where $\widetilde{W}$ is another space-time white noise and 
	\begin{equation}\label{sigmat}
		\widehat{\sigma}(x)= \theta ^{-\frac{1}{2}} \sigma (x) \mbox{ for every } x\in \mathbb{R}. 
	\end{equation}
	In other words, $v_{\theta}$ is the mild solution to the SPDE
	\begin{equation*}
		\frac{\partial v_{\theta }}{\partial t}(t,x)= - (-\Delta) ^{ \frac{\alpha}{2}}+ \widehat{\sigma} (v_{\theta}(t,x)) \dot{\widetilde{W}} (t,x), \hskip0.4cm t\geq 0, x\in \mathbb{R},
	\end{equation*}
	with vanishing initial condition and with  $\widehat{\sigma}$ from (\ref{sigmat}).  From Theorems \ref{tt1} and \ref{tt2}, we deduce the following limit theorems in $ L ^{1} (\Omega)$
	\begin{equation}\label{24f-1}
		V_{N,x} (v_{\theta}):=	N ^{-\frac{1}{\alpha}}\sum_{i=0} ^{N-1}\left( v_{\theta }(t_{i+1},x)-v_{\theta}(t_{i}, x)\right) ^{2}\to_{N \to \infty} C_{0, \alpha } ^{2} \int_{0}^{1} \widehat{\sigma }^{2} (v_{\theta }(s,x))ds,
	\end{equation}
	and if $ \frac{2\alpha}{\alpha -1}$ is an even integer,
	\begin{equation*}
		U_{N,x}(v_{\theta }):= \sum_{i=0} ^{N-1} \left| v_{\theta } (t_{i+1},x)-v_{\theta}(t_{i}, x)\right| ^{\frac{2\alpha}{\alpha -1}} \to _{N \to \infty}B_{0, \alpha} \int_{0} ^{1} \widehat{\sigma} ^{\frac{2\alpha}{\alpha-1}} (v_{\theta } (s,x))ds,
	\end{equation*}
	with $ C_{0, \alpha}, B_{0, \alpha}$  from  (\ref{deco}) and (\ref{boa}), respectively.  From these facts, we can easily obtain the behavior of the variations of $u_{\theta}$, the solution to the parametrized heat equation (\ref{1}).

	\begin{prop}\label{pp5}
		Let $ u_{\theta }$ be the solution to the SPDE (\ref{1}). Then
		\begin{equation}\label{10f-1}
			V_{N,x} (u_{\theta }):=N ^{-\frac{1}{\alpha}}\sum_{i=0} ^{N-1}\left( u_{\theta}(t_{i+1},x)-u_{\theta }(t_{i}, x)\right) ^{2}\to_{N \to \infty} C_{0, \alpha } ^{2} \theta ^{-\frac{1}{\alpha} }\int_{0}^{1} \sigma ^{2} (u_{\theta }(s,x))ds
		\end{equation}	
		and, if $ \frac{2\alpha}{\alpha -1} $ is an even integer,
		\begin{equation}\label{10f-2}
			U_{N,x} (u_{\theta }):= \sum_{i=0} ^{N-1} \left| u_{\theta } (t_{i+1},x)-u_{\theta}(t_{i}, x)\right| ^{\frac{2\alpha}{\alpha -1}} \to _{N \to \infty}B_{0, \alpha}\theta ^{-\frac{1}{\alpha -1}} \int_{0}^{1} \sigma ^{\frac{2\alpha}{\alpha -1}}(u_{\theta }(s,x)) ds.
		\end{equation}
		
	\end{prop}
	{\bf Proof: } Since $u_{\theta } (t,x)= v_{\theta } (\theta t, x)$, we have
	\begin{equation*}
		U_{N,x} (u_{\theta })= \sum_{i=0} ^{N-1} \left| v_{\theta } (\theta t_{i+1},x)-v_{\theta }(\theta t_{i}, x)\right| ^{\frac{2\alpha}{\alpha -1}}, 
	\end{equation*}
	and since $\theta t_{i}, i=1,...,N$ form a partition of the interval $[0, \theta]$, we get from (\ref{13f-3}),
	\begin{eqnarray*}
		&&	U_{N,x} (u_{\theta }) \to _{N \to \infty}B_{0, \alpha} \int_{0} ^{\theta} \widehat{\sigma} ^{\frac{2\alpha}{\alpha -1}}(v_{\theta }(s,x))ds \\
		&&=B_{0, \alpha}\theta ^{-\frac{\alpha}{\alpha -1}} \int_{0} ^{\theta} \sigma ^{\frac{2\alpha}{\alpha -1}}(v_{\theta }(s,x))ds =\theta ^{-\frac{\alpha}{\alpha -1}+1}\int_{0} ^{1} \sigma ^{\frac{2\alpha}{\alpha -1}}(v_{\theta }(\theta s,x))ds\\
		&&=B_{0, \alpha} \theta ^{-\frac{1}{\alpha -1}}\int_{0} ^{1} \sigma ^{\frac{2\alpha}{\alpha -1}}(u_{\theta }(s,x))ds.
	\end{eqnarray*} 
	Similarly, via (\ref{13f-2}),
	
	\begin{eqnarray*}
		V_{N, x} (u_{\theta}) &=& N ^{-\frac{1}{\alpha}}\sum_{i=0} ^{N-1} \left( v_{\theta}(\theta t_{i+1}, x) - v_{\theta}(\theta t_{i}, x) \right) ^{2} \\
		&\to _{N \to \infty} &C_{0, \alpha  }^{2} \theta ^{-\frac{1}{\alpha}}\int_{0} ^{\theta}\widehat{\sigma } ^{2} (v_{\theta}(s,x))ds = \theta ^{-\frac{1}{\alpha}} \theta ^{-1} \int_{0} ^{\theta}\sigma ^{2} (v_{\theta}(s,x))ds\\
		&=& C_{0, \alpha } ^{2}\theta ^{-\frac{1}{\alpha}} \int_{0}^{1} \sigma ^{2} (u_{\theta}(s,x))ds.
	\end{eqnarray*}
	\qed 
	
	\subsection{Estimation of the anomality parameter $\alpha$}\label{sec42}
	
	We first  illustrate how the parameter $\alpha$ of the fractional Laplacian (we will call it "the anomality parameter") can be estimated in the model (\ref{1}) via the result in Proposition \ref{pp5}.  We assume that we dispose on the observations $ \left( u_{\theta}(t_{i},x), i=0,1,...,N\right)$, where $x\in \mathbb{R}$ is arbitrary and $t_{i}, i=0,1,...,N$ are given by  (\ref{ti}). The estimation of this parameter, although it represents a natural question, seems to be new in the literature. 
	
	We will define the estimator 
	\begin{equation}
		\label{est3}
		\widehat{\alpha}_{N} = \frac{ \log(N)}{\log (A_{N})}, \hskip0.5cm N\geq 1
	\end{equation}
	with 
	\begin{equation}\label{an}
		A_{N}:= \frac{ \sum_{i=0} ^{N-1} \left( u_{\theta}( t_{i+1}, x)- u_{\theta}(t_{i}, x)\right) ^{2}}{ \frac{1}{N} \sum_{i=1} ^{N} \sigma ^{2} (u_{\theta}(t_{i}, x))}, \hskip0.4cm N\geq 1.
	\end{equation}
	
	We show below that (\ref{est3}) constitutes a consistent estimator for the anomality parameter $\alpha$. 
	\begin{prop}
		Let $ \widehat{\alpha}_{N}$ be given by (\ref{est3}). Then 
		\begin{equation*}
			\widehat{\alpha}_{N} \to_{N \to \infty} \alpha \mbox{ in probability. } 
		\end{equation*}
	\end{prop}
	{\bf Proof: } By the limit  (\ref{10f-1}) in Proposition \ref{pp5}, we have
	\begin{equation*}
		N ^{-\frac{1}{ \alpha}} A_{N} \to _{N \to \infty} C_{0, \alpha} ^{2} \theta ^{-\frac{1}{\alpha}}\mbox{ in probability, }
	\end{equation*}
	with $ C_{0, \alpha} $ from (\ref{deco}) and $A_{N}$ given by (\ref{an}). This implies that 
	\begin{equation*}
		E_{N}:= -\frac{1}{ \alpha}\log(N) + \log (A_{N})-2\log (C_{0, \alpha})+\frac{1}{\alpha }\log(\theta)\to _{N \to \infty }0 \mbox{ in probability. } 
	\end{equation*}
	Therefore
	\begin{equation*}
		\frac{ E_{N}}{\log (N)}=-\frac{1}{ \alpha} + \frac{ \log (A_{N})}{\log (N)} - \frac{ 2  \log (C_{0, \alpha})}{\log (N)}+ \frac{\log(\theta)}{\alpha \log(N)}\to _{N \to \infty}0, 
	\end{equation*}
	in probability. The last relation clearly gives 
	\begin{equation*}
		\frac{ \log (A_{N})}{\log (N)} \to _{N \to \infty} \frac{1}{ \alpha} \mbox{ in probability, }
	\end{equation*}
	and the conclusion follows. \qed 
	
	Notice that the estimator (\ref{est3}) does not depend on the other parameter of the model, the drift $\theta$. 
	
	\subsection{Estimation of the drift  parameter}

	The purpose is to identify the parameter $\theta$ in the SPDE (\ref{1}), by assuming that $\alpha$ is known and that, as above, we observe the solution $u_{\theta}$ at discrete time and at any arbitrary spatial point.

	Proposition \ref{pp5} leads naturally to the definition of two estimators for the parameter $  \theta $ Let us set, for $N\geq 1$, 
	\begin{equation}
		\label{est1}
		\widehat{\theta}_{N, 1}=\left( \frac{ V_{N, x} (u_{\theta}) }{C_{0, \alpha }^{2} \frac{1}{N}\sum_{i=1} ^{N} \sigma ^{2} (u(t_{i}, x))}\right) ^{-\alpha}
	\end{equation}
	and
	\begin{equation}
		\label{est2}
		\widehat{\theta}_{N,2 } =\left( \frac{ U_{N,x}(u_{\theta })}{B_{0, \alpha}\frac{1}{N}\sum_{i=1} ^{N} \sigma ^{\frac{2\alpha}{\alpha -1}} (u_{\theta }(t_{i}, x))}\right) ^{-(\alpha -1)} . 
	\end{equation}
	
	We notice again that the above estimators can be effectively computed if one dispose on the the observations $ (u_{\theta}(t_{i}, x), i=0, 1,...,N)$ with $t_{i}$ given by (\ref{ti}) and with $x\in \mathbb{R}$ arbitrary. We show below the consistency of the above estimators. 
	
	\begin{prop}
		Let $\widehat{\theta}_{N,1}$ and $\widehat{\theta}_{N,2}$ be given by (\ref{est1}), \ref{est2}), respectively. Then, for $j=1,2$,
		\begin{equation*}
			\widehat{\theta}_{N,j }\to_{N \to \infty} \theta \mbox{ in probability. }
		\end{equation*}
	\end{prop}
	{\bf Proof: }We start with the estimator $\widehat{\theta}_{N,1}$. We write 
	
	\begin{eqnarray*}
		\widehat{\theta}_{N,1} ^{-\frac{1}{\alpha}}-\theta ^{-\frac{1}{\alpha}}
		&=& \frac{ V_{N, x} (u_{\theta})}{ C_{0, \alpha } ^{2} \frac{1}{N}\sum_{i=1} ^{N} \sigma ^{2} (u_{\theta}(t_{i}, x))}- \frac{ V_{N, x}(u_{\theta})}{C _{0, \alpha } ^{2} \int_{0} ^{1} \sigma ^{2}(u_{\theta}(s,y))ds }\\
		&&+\frac{ V_{N, x}(u_{\theta})}{C _{0, \alpha } ^{2} \int_{0} ^{1} \sigma ^{2}(u_{\theta}(s,y))ds }- \theta ^{-\frac{1}{\alpha}}\\
		&=& \frac{ V_{N, x}(u_{\theta})}{C_{0, \alpha } ^{2}}\left( \frac{1}{\frac{1}{N}\sum_{i=1} ^{N} \sigma ^{2} (u_{\theta}(t_{i}, x)) }-\frac{1}{ \int_{0} ^{1} \sigma ^{2}(u_{\theta}(s,y))ds }\right)\\
		&&+ \frac{ V_{N, x}(u_{\theta})}{C _{0, \alpha } ^{2} \int_{0} ^{1} \sigma ^{2}(u_{\theta}(s,y))ds }- \theta ^{-\frac{1}{\alpha}}
	\end{eqnarray*}

	The first summand goes to zero in probability since the Riemann sum \\$ \sum_{i=1} ^{N} \sigma ^{2} (u_{\theta}(t_{i}, x))$ converges almost surely, as $ N\to \infty$, to the Lebesgue integral\\ $\int_{0} ^{1} \sigma ^{2} ( u_{\theta }(s,x))ds$ and $ V_{N, x}(u_{\theta}) $ converges in $ L^{1}(\Omega)$ to a nontrivial limit by Proposition \ref{pp5}.  The second  summand converges to $0 $ in probability as $N\to \infty$, due to (\ref{10f-1}).  Thus, the conclusion is obtained for $j=1$. The consistency of (\ref{est2}) follows in a similar way. \qed

	\begin{remark}\label{rem2}
		\begin{itemize}
			\item 	To estimate the drift  $\theta$ in (\ref{1}), we can use  both   sequences $V_{N,x}(u_{\theta})$ and $U_{N,x}(u_{\theta})$ and the limit behavior  in Theorems \ref{tt1} and \ref{tt2}.  The spatial version of the sequence $ U_{N,x}(u_{\theta})$ has been used in  \cite{GT} to identify the drift parameter. The sequence $ V_{N, x}(u_{\theta})$ provides an alternative estimator for the drift parameter but also it allows to identify the anomality parameter $\alpha$, see Section \ref{sec42}.
			
			\item As in \cite{BT1}, \cite{BT2} or \cite{GT}, one can propose a consistent  estimator based on the average of the temporal quadratic variations at different spatial points of the solution $u_{\theta}$.  An analysis of such an averaged estimator (based on a finite or infinite number of spatial points) and a comparaison with the estimators  (\ref{est1}), (\ref{est2}) would be of interest.

			\item If $\sigma \equiv 1$, then a Central Limit Theorem has been obtained for $\widehat{\theta}_{N,2}$ (see \cite{Mahtu2} or \cite{T2}).  The obtention of the  asymptotic distribution for the power variation in the nonlinear case constitutes an interesting problem, as well as the limit distribution of the estimators $\widehat{\alpha}_{N}$ and $\widehat{\theta}_{N, 1}$ (even in the linear noise case).
		\end{itemize}
	\end{remark}

	\section{Appendix}
	In this section, we included a description of the main properties of the Green kernel associated to the fractional Laplacian operator and the proof of a technical result.
	
	\subsection{Properties of the Green kernel}
	
	For the definition and other properties of the fractional Laplacian, one may consult, among others, the reference  \cite{Ja1}. Here, we will mainly  work with the Green kernel associated to this operator. The best way to define the Green kernel $G_{\alpha}$ is via its Fourier transform (we use $\mathcal{F}f$ to denote the Fourier transform of the  function $f$), namely 
	\begin{equation}\label{ga}
		\mathcal{F} G_{\alpha} (t, \cdot)(\xi):=\int_{\mathbb{R}} e ^{i x \xi }G_{\alpha} (t,x)dx= e ^{-t\vert \xi \vert ^{\alpha} }
	\end{equation}
	for $t\in [0, T] $ and $ \xi \in \mathbb{R}$, i.e.
	\begin{equation*}
		G_{\alpha} (t,x) =(2\pi ) ^{-1} \int_{\mathbb{R}} e ^{-ix\xi } e ^{-t\vert \xi \vert ^{\alpha}}d\xi.
	\end{equation*}
	
	Let us list some useful properties of the Green kernel (\ref{ga}).
	
	\begin{itemize}
		\item For every $t>0$, $G_{\alpha }(t, \cdot) $ is the density of the stable L\'evy process. In particular, this implies 
		\begin{equation*}
			\int_{\mathbb{R}} G_{\alpha }(t,x)dx =1. 
		\end{equation*}
		\item  For every $t>0$,  the kernel $ G_{\alpha}(t,x) $ is positive, real -valued and symmetric with respect to the space variable $x\in \mathbb{R}$.
		
		\item We have the scaling property
		\begin{equation*}
			G_{\alpha} (t,x) = t ^{-\frac{1}{\alpha}}G_{\alpha} (1, t ^{-\frac{1}{\alpha}} x)
		\end{equation*}
		for $t>0, x\in \mathbb{R}$. 
		
		\item There exists two constants $0< K_{\alpha } '< K_{\alpha}$ such that for every $t\in [0, T], x\in \mathbb{R}$, 
		\begin{equation}
			\label{3i-2}
			K_{\alpha } ' \frac{ t ^{-\frac{1}{\alpha}}}{ \left( 1+\vert t ^{ -\frac{1}{\alpha}}x\vert \right)^{1+\alpha}}\leq G_{\alpha } (t,x) \leq K_{\alpha} \frac{ t ^{-\frac{1}{\alpha}}}{ \left( 1+\vert t ^{ -\frac{1}{\alpha}}x\vert \right)^{1+\alpha}}.
		\end{equation}
	\end{itemize}

	\subsection{Proof of Proposition \ref{pp4}}
	We have
	\begin{eqnarray*}
		&&	\left( u(t+\delta, x)- u(t,x)\right) - \sigma (u( t(\delta), x))( \widetilde{\Delta} _{1} u_{0} ) (t, \delta) \\
		&=&\int_{ t(\delta)} ^{t+\delta } \int_{\mathbb{R}} \left( G_{\alpha} (t+\delta -a, x-y) - G_{\alpha} (t-a, x-y)\right)\\
		&& \left( \sigma (u(a,y)) - \sigma (u( t(\delta), x))\right) W (da, dy)\\
		&&+ \int_{0} ^{t(\delta)} \int_{\mathbb{R}}\left( G_{\alpha} (t+\delta -a, x-y) - G_{\alpha} (t-a, x-y)\right) \sigma (u(a,y)) W (da, dy)\\
		&&-  \int_{0} ^{t(\delta)} \int_{\mathbb{R}}\left( G_{\alpha} (t+\delta -a, x-y) - G_{\alpha} (t-a, x-y)\right) \sigma (u( t(\delta), x))\widetilde{W} (da, dy)\\
		&=& I_{1} + I_{2}+ I_{3}.
	\end{eqnarray*}
	Thus
	\begin{equation}\label{3f-21}
		\mathbf{E} \left| 	\left( u(t+\delta, x)- u(t,x)\right) - \sigma (u( t(\delta), x))( \widetilde{\Delta} _{1} u_{0} ) (t, \delta) \right| ^{2}\leq C( 	\mathbf{E} I_{1} ^{2} + 	\mathbf{E} I_{2} ^{2} + 	\mathbf{E} I_{3} ^{2}).
	\end{equation}
	
	For the term $I_{2}$,  by the isometry of the Dalang-Walsh integral (\ref{iso-dw}),
	\begin{eqnarray}
		\mathbf{E} I_{2} ^{2} &=& \int_{0} ^{t(\delta)}da  \int_{\mathbb{R}}dy \left( G_{\alpha} (t+\delta -a, x-y) - G_{\alpha} (t-a, x-y)\right)^{2} 	\mathbf{E} \sigma ^{2}  (u(a,y)) \nonumber\\
		&\leq & C_{T}   \int_{0} ^{t(\delta)}da  \int_{\mathbb{R}}dy \left( G_{\alpha} (t+\delta -a, x-y) - G_{\alpha} (t-a, x-y)\right)^{2} \nonumber \\
		&\leq & C_{T} \delta ^{2- \frac {\beta(\alpha+1)}{\alpha}}=C_{T} \delta ^{\frac{2\alpha}{2\alpha+1}},\label{3f-22}
	\end{eqnarray}
	where we used succesively the inequality (\ref{3i-5}) and Lemma \ref{ll4}. In a entirely similar way 
	\begin{equation}\label{3f-23}
		\mathbf{E} I_{3} ^{2} \leq C_{T}  \delta ^{\frac{2\alpha}{2\alpha+1}}.
	\end{equation}
	We will estimate the term $ I_{1}$ as follows. We use the isometry property (\ref{iso-dw}), the Lipschitz assumption on $\sigma $ and (\ref{2i-2}) to obtain
	\begin{eqnarray*}
		\mathbf{E} I_{1} ^{2} &=& \int_{ t(\delta)} ^{t+\delta }da  \int_{\mathbb{R}}dy  \left( G_{\alpha} (t+\delta -a, x-y) - G_{\alpha} (t-a, x-y)\right) ^{2} \\
		&&	\mathbf{E} \left( \sigma (u(a,y)) - \sigma (u( t(\delta), x))\right) ^{2} \\
		&\leq &  C_{T} \int_{ t(\delta)} ^{t+\delta }da  \int_{\mathbb{R}}dy  \left( G_{\alpha} (t+\delta -a, x-y) - G_{\alpha} (t-a, x-y)\right) ^{2}\\
		&& \left[   \vert a-t(\delta)\vert ^{ 1-\frac{1}{\alpha}} + \vert x-y\vert ^{\alpha -1}   \right]\\
		&\leq & C_{T} \int_{ t(\delta)} ^{t+\delta }da  \int_{\mathbb{R}}dy  \left( G_{\alpha} (t+\delta -a, x-y) - G_{\alpha} (t-a, x-y)\right) ^{2} \\
		&&\left[ \vert \delta +\delta ^{\beta } \vert ^{1-\frac{1}{\alpha}} + \vert x-y\vert ^{\alpha -1} \right] \\
		&\leq & C_{T} \int_{ t(\delta)} ^{t+\delta }da  \int_{\mathbb{R}}dy  \left( G_{\alpha} (t+\delta -a, x-y) - G_{\alpha} (t-a, x-y)\right) ^{2}\\
		&&\left[  \delta ^{ \beta (1-\frac{1}{\alpha})}+ \vert x-y\vert ^{\alpha -1} \right],
	\end{eqnarray*}
	where we used $\beta <1$ for the last inequality. Next, since
	$$ \left( G_{\alpha} (t+\delta -a, x-y) - G_{\alpha} (t-a, x-y)\right) ^{2}\leq 2 \left( G^{2}_{\alpha} (t+\delta -a, x-y)+ G^{2} _{\alpha} (t -a, x-y)\right)$$
	and using $G_{\alpha } (t,x)=0$ if $t<0$,
	\begin{eqnarray*}
		\mathbf{E} I_{1} ^{2}&\leq & C_{T} \delta ^{\beta (1-\frac{1}{\alpha})} \int_{ t(\delta)} ^{t+\delta} da \int_{\mathbb{R}} dy G^{2}_{\alpha} (t+\delta -a, x-y)\\
		&&+ C_{T}  \delta ^{\beta (1-\frac{1}{\alpha})} \int_{ t(\delta)} ^{t} da \int_{\mathbb{R}} dy G^{2}_{\alpha} (t -a, x-y)\\
		&&+C_{T} \int_{ t(\delta)} ^{t+\delta} da \int_{\mathbb{R}} dy  G^{2}_{\alpha} (t+\delta -a, x-y)\vert x-y\vert ^{\alpha -1} \\
		&&+ C_{T} \int_{ t(\delta)} ^{t} da \int_{\mathbb{R}} dy G_{\alpha} ^{2} (t-a, x-y)\vert x-y\vert ^{\alpha -1} .
	\end{eqnarray*}
	By Lemma 2 in \cite{GT}, we have, for $0\leq s\leq t$, 
	\begin{equation}
		\label{gt1}
		\int_{s} ^{t} da\int_{\mathbb{R}}dy G_{\alpha} ^{2} (t-a, y)\leq C (t-s) ^{1-\frac{1}{\alpha}}.
	\end{equation}
	By using (\ref{gt1}),
	\begin{eqnarray*}
		\mathbf{E} I_{1} ^{2}&\leq & C_{T}  \delta ^{2\beta (1-\frac{1}{\alpha})}\\
		&&+C_{T} \int_{ t(\delta)} ^{t+\delta} da \int_{\mathbb{R}} dy  G^{2}_{\alpha} (t+\delta -a, x-y)\vert x-y\vert ^{\alpha -1} \\
		&&+ C_{T} \int_{ t(\delta)} ^{t} da \int_{\mathbb{R}} dy G_{\alpha} ^{2} (t-a, x-y)\vert x-y\vert ^{\alpha -1} .
	\end{eqnarray*}
	The last summand has been already treated in the proof of Proposition 4 in \cite{GT}.  By these estimates
	$$\int_{ t(\delta)} ^{t} da \int_{\mathbb{R}} dy G_{\alpha} ^{2} (t-a, x-y)\vert x-y\vert ^{\alpha -1} \leq C  \delta ^{2\beta (1-\frac{1}{\alpha})}$$
	and analogously,
	$$ \int_{ t(\delta)} ^{t+\delta} da \int_{\mathbb{R}} dy  G^{2}_{\alpha} (t+\delta -a, x-y)\vert x-y\vert ^{\alpha -1} \leq C  \delta ^{2\beta (1-\frac{1}{\alpha})}.$$
	Thus,
	\begin{equation}\label{3f-24}
		\mathbf{E} I_{1}^{2} \leq C  \delta ^{2\beta (1-\frac{1}{\alpha)}}= \delta ^{\frac{4(\alpha-1)}{2\alpha+1}}.
	\end{equation}
	The conclusion is obtained by plugging the estimates (\ref{3f-22}), (\ref{3f-23}) and (\ref{3f-24}) into (\ref{3f-21}).  \qed
	
	The following lemma has been used in the proof of Proposition \ref{pp4}. 
	
	\begin{lemma}\label{ll4}
		For every $0\leq s\leq t\leq T$, 
		\begin{equation*}
			I(s,t, \delta)=	\int_{0} ^{s}da \int_{\mathbb{R}}dy \left( G_{\alpha} (t+\delta -a, x-y) - G_{\alpha} (t-a, x-y)\right) ^{2} \leq C _{T} \delta ^{2} (t-s)^{-\frac{\alpha+1}{\alpha}}. 
		\end{equation*}
	\end{lemma}
	{\bf Proof: } By using Parseval's identity,
	\begin{eqnarray*}
		I(s, t, \delta)&=& \int_{0} ^{s} da \int_{\mathbb{R}} d\xi \left| e ^{-(t+\delta -a) \vert \xi \vert ^{\alpha}}-e ^{-(t-a) \vert \xi \vert ^{\alpha}}\right| ^{2} \\
		&=&  \int_{0} ^{s} da \int_{\mathbb{R}} d\xi e ^{-2(t-a) \vert \xi \vert ^{\alpha}}\left| 1- e ^{-\delta \vert \xi \vert ^{\alpha} }\right| ^{2}.
	\end{eqnarray*}
	Next, since$\left| 1- e ^{-2s\vert \xi \vert ^{\alpha}}\right|\leq 1$, 
	\begin{eqnarray*}
		I(s, t,\delta) &\leq & \delta ^{2} \int_{0} ^{s} da \int_{\mathbb{R}} d\xi e ^{-2(t-a)\vert \xi \vert ^{\alpha}}\vert \xi \vert ^{2\alpha} \\
		&=& C \delta ^{2} \int_{\mathbb{R}} d\xi \vert \xi \vert ^{\alpha} e ^{-2(t-s) \vert \xi \vert ^{\alpha}}\left( 1- e ^{-2s\vert \xi \vert ^{\alpha}}\right)\\
		&\leq & C \delta ^{2} \int_{\mathbb{R}} d\xi \vert \xi \vert ^{\alpha}e ^{-2(t-s) \vert \xi \vert ^{\alpha}}\\\
		&\leq & C \delta ^{2} (t-s) ^{-\frac{\alpha+1}{\alpha}}
	\end{eqnarray*}
	where we performed the change of variables $ (t-s) ^{\frac{1}{\alpha} } \xi = \tilde{\xi}$. \qed


\begin{thebibliography}{99}
		
		\bibitem{ACP}
		{R. Altmeyer. I.  Cialenco and G. Pasemann (2023): }
		{\em Parameter estimation for semilinear SPDEs from local measurements. } Bernoulli 29 (3),  2035–2061.
		
		
		\bibitem{AT1}
		{O. Assaad and C. A. Tudor (2021): }{\em Pathwise analysis and parameter estimation for the stochastic Burgers equation.}  Bull. Sci. Math. 170, Paper No. 102995, 26 pp. 
		
		
		
		\bibitem{ANTV}
		{O. Assaad, D. Nualart, C.A. Tudor and L. Viitasaari (2021): }	{\em Quantitative normal approximations for the stochastic fractional heat equation. }Stoch. Partial Differ. Equ. Anal. Comput. 10 (2022), no. 1, 223–254.
		
		
		\bibitem{BEV}
		S. Bajja, K. Es-Sebaiy,  and L. Viitasaari (2018): {\em Limit theorems for quadratic variations of the Lei-Nualart process.} Stochastic processes and applications, 105–121, Springer Proc. Math. Stat., 271, Springer, Cham.
		
		\bibitem{BT1}
		M. Bibinger, M.  Trabs (2020):  {\em Volatility estimation for stochastic PDEs using high-frequency observations. }
		Stochastic Process. Appl., 130 (5), 3005-3052.
		
		
		\bibitem{BT2}
		{ M. Bibinger and M. Trabs (2019): }{\em  On central limit theorems for power variations of the solution to the stochastic heat equation. }Stochastic models, statistics and their applications, 69–84, Springer Proc. Math. Stat., 294, Springer, Cham.
		
		
		
		\bibitem{Chong}
		{C. Chong (2020): }{\em High-frequency analysis of parabolic stochastic PDEs. } Ann. Statist., 48 1143-1167.
		
		
		\bibitem{CD}
		{C. Chong and R. Dalang (2023): }{\em  Power variations in fractional Sobolev spaces for a class of parabolic stochastic PDEs. } Bernoulli 29 (3), 1792-1820.
		
		
		
		\bibitem{Cia}
		{I. Cialenco (2018): }{\em Statistical inference for spdes: an overview. }Statistical Inference for Stochastic Processes, 21(2), 309-329.
		
		
		
		
		
		
		
		\bibitem{CDK}
		{I. Cialenco, F.  Delgado-Vences, H.-J.  Kim (2020): }
		{\em Drift estimation for discretely sampled SPDEs. } Stoch. PDE: Anal. Comp., 8 (2020),  895-920.
		
		\bibitem{CH}
		{I. Cialenco and Y. Huang (2020): }{\em  A note on parameter estimation for discretely sampled spdes. }Stoch. Dyn. 20 (2020), no. 3, 2050016, 28 pp.
		
		\bibitem{CK}
		{	I. Cialenco, H-J. Kim  (2022): }
		{\em 	Parameter estimation for discretely sampled stochastic heat equation driven by space-only noise. }  Stochastic Process. Appl. 143, 1–30.
		
		
		\bibitem{CKP}
		{I. Cialenco, H.-J. Kim,  and G. Pasemann  (2024): }{\em  Statistical analysis
			of discretely sampled semilinear SPDEs: a power variation approach. } Stoch. Partial Differ. Equ. Anal. Comput. 12 (1), 326-351.
		
		
		
		\bibitem{Da}
		{R. Dalang  (1999: }{\em Extending the Martingale Measure Stochastic Integral With Applications to Spatially Homogeneous S.P.D.E.’s. }Electron. J. Probab. Volume 4, paper no. 6, 29 pp.
		
		
		
		\bibitem{DD}
		L. Debbi and M. Dozzi (2005): {\em On the solutions of nonlinear stochastic fractional partial differential
			equations in one spatial dimension. }Stoch. Proc. Appl., 115: 1761-1781.
		
		
		\bibitem{GT}
		{J. Gamain and C. A. Tudor (2023): }{\em Exact variation and drift parameter estimation for the nonlinear fractional stochastic heat equation. } Jpn. J. Stat. Data Sci. 6 (1), 381-406.
		
		
		
		
		
		
		
		\bibitem{G1}
		S. Gaudlitz  (2023): {\em Non-parametric estimation of the reaction term in semi-linear SPDEs with spatial ergodicity. } Preprint arXiv:2307.05457.
		
		\bibitem{GR}
		S. Gaudlitz, M. Reiss (2023): {\em  Estimation for the reaction term in semi-linear SPDEs under small diffusivity. } Bernoulli, 29 (4), 3033 - 3058.
		
		\bibitem{HiTr}
		{F.  Hildebrandt  and M. Trabs  (2021): }{\em Parameter estimation for SPDEs based on discrete observations in time and space. }Electron. J. Stat. 15 (1), 2716-2776. 
		
		
		
		
		\bibitem{Ja1}
		{N. Jacob and H.G. Leopold (1993): }{\em  Pseudo differential operators with variable order of differentiation generating Feller semigroups. }  Integral Equations Operator Theory,  17, 544-553.
		
		\bibitem{KU}
		{Y. Kaino and M.  Uchida (2021): }
		{\em Parametric estimation for a parabolic linear SPDE model based on sampled data. }
		J. Statist. Plann. Inference, 211, 190-220.
		
		
		
		
		
		
		
		
		\bibitem{MahTu}
		{Z. Khalil-Mahdi and C. A. Tudor (2019): }{ \em On the distribution and $q$-variation of the solution to the heat equation with Fractional Laplacian. } Probability Theory and Mathematical Statistics, 39 (2), 315–335.
		
		
		\bibitem{Mahtu2}
		Z. Mahdi Khalil and C. A. Tudor (2019). Estimation of the drift parameter for
		the fractional stochastic heat equation via power variation, Mod. Stoch.
		Theory Appl. 6 (4), 397-417.
		
		\bibitem{NNT}
		{	I. Nourdin, D. Nualart and C. A. Tudor (2010): }
		{\em 	Central and non-central limit theorems for weighted power variations of fractional Brownian motion. }
		Ann. Inst. Henri Poincaré Probab. Stat.   46 (4), 1055–1079.
		
		
		\bibitem{NP-book}
		{I. Nourdin and G. Peccati (2012): }{ Normal Approximations with
			Malliavin Calculus From Stein's Method to Universality}. Cambridge
		University Press.
		
		
		
		\bibitem{PoTr}
		{J. Pospisil and R. Tribe (2007): }{\em Parameter estimates and exact variations for stochastic heat equation heat equations driven by space-time white noise.  }Stochastic  Analysis  and Applications 25(3), 593-611.
		
		
		\bibitem{Sw}
		J. Swanson (2007): {\em  Variations of the solution to a stochastic
			heat equation. } The Annals of Probability,  35(6), 2122-2159.
		
		\bibitem{TKU}
		Y. Tonaki, Y. Kaino and  M. Uchida  (2023): {\em  Parameter estimation for a linear parabolic SPDE model in two space dimensions with a small noise. }Stat. Inference Stoch. Process. 27 (1), 123-179.
		
		\bibitem{T2}
		C. A. Tudor (2022): Stochastic Partial Differential Equations with Additive Gaussain noise: Analysis and Inference. World Scientific. 
		
		
		\bibitem{Walsh} 
		J. W  Walsh  (1986): {\em An introduction to stochastic partial differential equations.}  
		\'{E}cole d' \'{e}t\'{e} de probabilit\'{e}s de Saint-Flour, XIV—1984,  265–439, LNM 1180, Springer.
		
		\bibitem{ZZ}
		{M. Zili and E. Zougar (2019): }{\em Exact variations for stochastic heat equations with piecewise constant coefficients and applications to parameter estimation. } Teor. Ĭmovīr. Mat. Stat. No. 100 (2019), 75–101.
		
		
	\end{thebibliography}
\end{document}